\documentclass[reqno,10.5 pt]{amsart}
\usepackage{amsmath,amssymb,latexsym}
\usepackage[colorlinks=true, linkcolor={black}]{hyperref}
\usepackage{amsthm}
\usepackage{comment}
\usepackage{color}  
\usepackage[all]{xy}
\usepackage{tikz-cd}
\usepackage{multirow}
\usepackage{multicol}
\usepackage{longtable}
\usepackage{caption}
\usetikzlibrary{cd}
\hypersetup{
    colorlinks=true, 
    linktoc=all,     
    linkcolor=blue,  
    urlcolor=black,
}

\newcommand{\F}{{\mathbb{F}}}
\newcommand{\Z}{{\mathbb{Z}}}
\newcommand{\Q}{{\mathbb{Q}}} 
\newcommand{\QQ}{\overline{{\mathbb{Q}}}}
\newcommand{\N}{{\mathbb{N}}}    
\newcommand{\R}{{\mathbb{R}}}    

\newcommand{\q}{{\mathfrak{q}}}
\newcommand{\p}{{\mathfrak{p}}}

\newcommand{\OO}{{\mathcal{O}}}

\newcommand{\Ker}{\mathrm{Ker}}

\makeatletter
\newcommand*{\rom}[1]{\expandafter\@slowromancap\romannumeral #1@}
\makeatother
\newcommand\isomto{\stackrel{\sim}{\smash{\longrightarrow}\rule{0pt}{0.4ex}}}
\DeclareFontFamily{U}{wncy}{}
    \DeclareFontShape{U}{wncy}{m}{n}{<->wncyr10}{}
    \DeclareSymbolFont{mcy}{U}{wncy}{m}{n}
    \DeclareMathSymbol{\Sh}{\mathord}{mcy}{"58}
		
\theoremstyle{plain}

\newtheorem{theorem}{Theorem}[section]
\newtheorem*{theorem*}{Theorem}

\newtheorem{proposition}[theorem]{Proposition}
\newtheorem{rem}[theorem]{Remark}
\newtheorem{lemma}[theorem]{Lemma}
\newtheorem{corollary}[theorem]{Corollary}
\newtheorem{defn}[theorem]{Definition}

\begin{document}
\title{$\sqrt{-3}$-Selmer groups,  ideal class groups and large $3$-Selmer ranks} 

\author{Somnath Jha, Dipramit Majumdar and Pratiksha Shingavekar}
\address{\begin{scriptsize}Somnath Jha, \href{mailto:jhasom@iitk.ac.in}{jhasom@iitk.ac.in}, Indian Institute of Technology Kanpur, India\end{scriptsize}}
\address{\begin{scriptsize}Dipramit Majumdar, \href{mailto:dipramit@iitm.ac.in}{dipramit@iitm.ac.in}, Indian Institute of Technology Madras, India\end{scriptsize}}
\address{\begin{scriptsize}Pratiksha Shingavekar, \href{mailto:pshingavekar@gmail.com}{pshingavekar@gmail.com}, Chennai Mathematical Institute, India\end{scriptsize}}

\begin{abstract} 
  We consider the family of elliptic curves $E_{a,b}:y^2=x^3+a(x-b)^2$ with $a,b \in \Z$. These elliptic curves have a rational $3$-isogeny, say $\varphi$. We give an upper and a lower bound on the rank of the $\varphi$-Selmer group of $E_{a,b}$ over $K:=\Q(\zeta_3)$ in terms of the $3$-part of the ideal class group of certain quadratic extension of $K$. Using our bounds on the Selmer groups, we construct infinitely many curves in this family with arbitrary large $3$-Selmer rank over $K$ and no non-trivial $K$-rational point of order $3$. We also show that for a positive proportion of natural numbers $n$, the curve $E_{n,n}/\Q$  has root number $-1$ and $3$-Selmer rank $=1$.
\end{abstract}

\maketitle
\section*{Introduction}\label{intro}
The $n$-Selmer group of an elliptic curve $E$ defined over a number field $F$ encapsulates various important arithmetic properties of $E$. Indeed, via the $n$-descent exact sequence 
\begin{equation}\label{eq:ndescentseq}
0 \longrightarrow {E(F)}/{nE(F)} \longrightarrow {\rm Sel}^{n}(E/F) \longrightarrow \Sh(E/F)[n] \longrightarrow 0,
\end{equation}
it is closely related to the Mordell-Weil and the Tate-Shafarevich groups of the elliptic curve.  A similar consideration holds, more generally, for an isogeny induced Selmer group of an elliptic curve (see \eqref{eq:defofsha}). Several important results on $n$-Selmer group of elliptic curves have been established recently by Bhargava et.al, Mazur-Rubin and many others.  Soon after Cassels defined the Selmer group of an elliptic curve, he considered \cite{cass3} the relation between a (degree $3$) isogeny induced Selmer group and the Tate-Shafarevich group of an elliptic curve over $\Q$ with  the $3$-part of the ideal class group of certain quadratic extension of $\Q$. There are many results  comparing the Selmer group of an elliptic curve with the ideal class group of a suitable number field including \cite{bk, ces, li, ss}.

For a number field $F$ and a finite set $S$ of its finite primes, we let $\OO_{S}$ be the ring of $S$-integers in $F$ and $Cl_{S}(F)$ be its $S$-ideal class group. We denote by $h^3_S$ the $3$-rank of the $Cl_{S}(F)$ i.e. $h^3_{S}:=\dim_{\F_3}Cl_{S}(F) \otimes \F_3 = \dim_{\F_3}Cl_{S}(F)[3]$. Throughout this paper, we set $K:=\Q(\zeta_3)$, where $\zeta_3=\frac{-1+\sqrt{-3}}{2}$ is a primitive cube root of unity. 

Now, let $E$ be an elliptic curve over $\Q$ with a rational isogeny $\varphi: E \to \widehat{E}$ of degree $3$. Then (up to a change of coordinates), the equation of $E$ can be written either in the form (i) $E_a: y^2= x^3+a$ with $a \in \Z$, $a \neq 0$ or (ii) $E_{a,b}:y^2=x^3+a(x-b)^2$ with $a, b \in \Z$, $ab(4a+27b) \neq 0$ (see \cite{jms} for details).  In \cite{jms}, we studied the $\varphi$-Selmer groups ${\rm Sel}^\varphi(E_a/K)$ of the curves $E_a: y^2=x^3+a$ over $K$.  
In fact, the $\varphi$-Selmer groups of $E_a$ have also been studied by Bandini, \cite{bes}, Nekov\'a\v r  etc., whereas the study of the $\varphi$-Selmer groups of the curves $E_{a,b}$ remains somewhat unexplored. In the first part of this article, we will complement the study  on $\varphi$-Selmer rank bounds of $E_a$  in \cite{jms} by exploring the $\varphi$- and $3$-Selmer groups of the curves $E_{a,b}:y^2=x^3+a(x-b)^2$ over $K$, denoted by ${\rm Sel}^{\varphi}(E_{a,b}/K)$ and ${\rm Sel}^{3}(E_{a,b}/K)$ respectively. The main goal of this article is to give an upper and a lower bound on the $\F_3$-ranks of ${\rm Sel}^{\varphi}(E_{a,b}/K)$, in terms of $h^3_{S}$ for certain $S$-class groups of $K(\sqrt a)$. Here the sets $S$ are certain subsets of the set of primes of bad reduction of $E_{a,b}$ over $K$ (Definition \ref{defofS1S2S3}). We then use these bounds for the Selmer groups induced by $\varphi$ and by it's dual isogeny $\widehat{\varphi}$ to get the upper and lower bounds for the $3$-Selmer group ${\rm Sel}^3(E_{a,b}/K)$.
 The main result of this article is as stated below (cf. Theorem \ref{cortotype2} for details).
\begin{theorem}\label{cortotype2i}
Let $E_{a,b}: y^2=x^3+a(x-b)^2$ with $a,b \in \Z$ and $ab(4a+27b) \neq 0$ and $\psi_{a,b}:E_{a,b} \to \widehat{E}_{a,b}$ be the rational isogeny of degree $3$, defined in \eqref{eq:defofpsiab}. Let $S_1$, $S_2$ and $S_3$ (as in the Definition \ref{defofS1S2S3}) be certain subsets of the set of primes of bad reduction of $E_{a,b}$ over $K$. Assume that $a \notin K^{*2}$ and put $L=L_a:=K[X]/(X^2-a)$. Denote by $S_{i,j}(L)$ the set of primes of $L$ lying above $S_i \cup S_j$. If $3 \nmid a$, then 
\begin{multline*}
\text{max } \big\{h^3_{S_{1,2}(L)}, \enspace h^3_{S_{1,3}(L)} +|S_3|-|S_2|-1 \big\} \leq \dim_{\F_3} {\rm Sel}^{\psi_{a,b}}(E_{a,b}/K) \\
\leq \text{min } \big\{h^3_{S_{1,2}(L)} +|S_{1,2}(L)|+|S_3|-|S_2|+1, \enspace h^3_{S_{1,3}(L)}+|S_{1,3}(L)|+2 \big\}.
\end{multline*}
\noindent Further, the global root number $\omega(E_{a,b}/K)$ is given by $\omega(E_{a,b}/K) = (-1)^{|S_2|+|S_3|+1}$.

In particular, if $S_1$, $S_2$ and $S_3$ are all empty, then $\dim_{\F_3} {\rm Sel}^{\psi_{a,b}}(E_{a,b}/K) \in \big\{h^3_L, \ h^3_L+1 \big\}$ and the global root number $\omega(E_{a,b}/K)=-1$. \qed
\end{theorem}
The bounds obtained here are indeed sharp and we have  explicit numerical examples (see Example (2), \S\ref{type2global}) where the upper and the lower bound differ by $1$.

We consider a subfamily  $E_{n,n}$, where $n \in \N$ and  make suitable choices of $n$, such that $\text{Sel}^{\psi_{n,n}}(E_{n,n}/\Q)$ vanishes. Using this, we prove in Proposition \ref{3.10} that for at least 6\% of square-free natural numbers $n$, the $3$-rank of the $3$-Selmer group of $E_{n,n}$ over $\Q$  i.e.  $\dim_{\F_3} \text{Sel}^{3}(E_{n,n}/\Q)$ is equal to $1$.

We now briefly explain the key ideas  behind the proof of Theorem \ref{cortotype2i}. The existence of a rational $3$-isogeny $\psi_{a,b}: E_{a,b} \rightarrow \widehat{E}_{a,b}$ implies that the residual representation $E_{a,b}(\QQ)[3]$ is a reducible $G_\Q$-module and gives rise to a $G_\Q$-stable order $3$ subgroup $(0, \pm \sqrt{a}b)$ which naturally corresponds to  the quadratic \'etale algebra $L_a:=K[X]/(X^2-a)$. The Galois cohomology group $H^1(G_K, E_{a,b}[\psi_{a,b}])$ can be identified with    the kernel of the norm maps $L^*_a/{L^{*3}_a} {\longrightarrow} K^*/{K^{*3}}$  (Prop. \ref{cohom}) and  one obtains an explicit description of the isogeny induced Selmer group ${\rm Sel}^{\psi_{a,b}}(E_{a,b}/K)$ in \eqref{eq:newseldef}.  We compute the images of the local Kummer maps and use some results in algebraic number theory  (see \S \ref{ant},  Props. \ref{boundsgen} and  \ref{newM1M2sel}) to first arrive  at some  preliminary bounds for ${\rm Sel}^{\psi_{a,b}}(E_{a,b}/K)$ (Theorem \ref{containment2}) via certain $S$-ideal class groups of $L_a$. Further, we obtain sharper bounds for ${\rm Sel}^{\psi_{a,b}}(E_{a,b}/K)$ in Theorem \ref{cortotype2}. This is achieved by using the formula of Cassels \eqref{eq:casselsfor}, that involves the study of local Tamagawa numbers to compare the sizes of $\psi_{a,b}$- and $\widehat{\psi}_{a,b}$-Selmer groups of $E_{a,b}$ and $\widehat{E}_{a,b}$, respectively.

We then use these upper and lower bounds on ${\rm Sel}^{\psi_{a,b}}(E_{a,b}/K)$ in Theorem \ref{cortotype2} to exhibit a sub-family of the curves $E_{a,b}$ that has arbitrary large $3$-Selmer ranks in Theorem \ref{large3selmer}. Producing families of elliptic curves with large Selmer rank has been of considerable interest in the literature. Cassels \cite{cass3} produced an infinite family of elliptic curves $E/\Q$ satisfying $j(E)=0$, $E(\Q)[3] \neq 0$ and such that $\dim_{\F_3} {\rm Sel}^3(E/\Q)$ is unbounded as $E$ varies. On the other hand, \cite{ks} have produced families of elliptic curves with large $n$-Selmer groups over $\Q$, for  $n \geq 5$. There is a recent important work of  Bhargava et. al. \cite{bkos} on the arbitrarily large  growth of $3$-torsion of the Tate-Shafarevich groups of elliptic curves over number fields  with cyclic $9$-isogenies. See Remark \ref{largecomp} for details on related works. We prove the following result in this direction.
\begin{theorem}
Given any $n \geq 0$, there exist infinitely many elliptic curves $E_{a,1}:y^2=x^3+a(x-1)^2$, with $a \in \Z$  such that $\dim_{\F_3} {\rm Sel}^3(E_{a,1}/K) \geq \dim_{\F_3} {\rm Sel}^{\psi_{a,1}}(E_{a,1}/K) \geq 2n$. These curves have a rational $3$-isogeny and no non-trivial $3$-torsion points over $K$.  
\end{theorem}

The article consists of $4$ sections and is structured as follows: In \S \ref{iso}, we state some preliminary definitions and recall some algebraic results from \cite{jms} which will be used later. In \S \ref{type2curves}, we will compute the images of local Kummer maps at each prime of $K$. Our main results on the $\F_3$-rank bounds on the Selmer groups (and the proof of Proposition \ref{3.10}) is contained in \S \ref{type2global}. The result on arbitrary large $3$-Selmer ranks is discussed in \S \ref{largeselmerrank}. Table \ref{tab:type2examples} contains examples of the curves $E_{a,b}$ for which we compute the bounds on their $\psi_{a,b}$- and $3$-Selmer groups over $K$ using our results from \S\ref{type2global}.

\section{Preliminaries}\label{iso}
\subsection{The curves $E_{a,b}$ and the $3$-isogenies}
The rational degree $3$-isogeny corresponds to a $G_\Q$-invariant subgroup $C$ of $E_{a,b}(\QQ)$ of order $3$ and given $C$,  the rational $3$-isogeny $\varphi$ is obtained as $\varphi: E_{a,b} \to E_{a,b}/C$. Via a change of variables $(x,y) \mapsto (\frac{x}{9}-\frac{4a}{3},\frac{y}{27})$, we can identify $E_{a,b}/C$ with $\widehat{E}_{a,b}:=E_{-27a,4a+27b}$. Thus, we obtain a rational $3$-isogeny $\psi_{a,b}:E_{a,b} \to \widehat{E}_{a,b}$ given explicitly using \cite{velu} as
\begin{equation}\label{eq:defofpsiab}
\psi_{a,b}(x,y) = \Big( \frac{9(x^3+\frac{4}{3}ax^2-4abx+4ab^2)}{x^2}, \frac{27y(x^3+4abx-8ab^2)}{x^3} \Big).
\end{equation}
Putting $d=4a+27b$, we obtain the corresponding dual isogeny $\widehat{\psi}_{a,b}: \widehat{E}_{a,b} \to E_{a,b}$ given by
\begin{equation}\label{eq:defofpsihatab}
\widehat{\psi}_{a,b}(x,y)  = \Big( \frac{x^3-36ax^2 +108adx -108ad^2}{3^4x^2}, \frac{y(x^3-108adx+216ad^2)}{3^6x^3} \Big).
\end{equation}
Set $\zeta_3:=\frac{-1+\sqrt{-3}}{2}$ and $K:=\Q(\zeta_3)$, and note that there is an isomorphism $\lambda_c: E_{ac^2,bc^2} \isomto E_{a,b}$ defined by $(x,y) \mapsto (c^{-2}x,c^{-3}y)$ over $K$ for any $c \in K^*$. Thus, without any loss of generality, we can assume that $gcd(a,b)$ is square-free in $K$. Further, observe that $\widehat{\psi}_{a,b}=\lambda_{27} \circ \psi_{-27a,d}$ and hence one has an isomorphism of Selmer groups ${\rm Sel}^{\widehat{\psi}_{a,b}}(\widehat{E}_{a,b}/K) \cong {\rm Sel}^{{\psi}_{-27a,d}}(E_{-27a,d}/K)$. 

\vspace{2mm}
\noindent {\bf{Notation}: } For any finite set $X$,  $|X|$ denotes its cardinality.
\begin{itemize}
\item  Given a multiplicative group $G$, define the subgroup $G^3:=\{x^3\mid x \in G\}$.
\item  For abelian groups $A$, $A'$ and a homomorphism $\eta: A \rightarrow A'$, let $A[\eta]$ denote the kernel of $\eta$. 
\item  $F$ denotes an arbitrary number field and let $\Sigma_F$ be the set of all its finite places. For a finite subset $S$ of $\Sigma_F$, $\OO_S$ denotes the ring of $S$-integers in $F$. For $\varpi \in \Sigma_F$, $F_\varpi$ will denote the completion of $F$ at $\varpi$ and its ring of integers will be denoted by $\OO_{F_\varpi}$.
\item  $\zeta_3=\frac{-1+\sqrt{-3}}{2}$ is the primitive third root of unity, $K=\Q(\zeta_3)$ and $\OO_K=\Z[\zeta_3]$ the ring of integers of $K$.  We will denote a finite prime of $\OO_K$ by $\q$.
\item $(\p)=(1-\zeta_3)$ is the unique prime of $K$ above $3$ with uniformizer $\p=1-\zeta_3$ in $\OO_{K_\p}$. By abuse of notation, we simply write $\p$ for the ideal $(\p)$ of $\OO_K$.
\item We set $L=L_a=\frac{K[X]}{(X^2-a)}$. Thus $L/K$ is a quadratic field extension  if $ a \notin K^{*2}$. Otherwise, $L \cong K \times K$ if $ a \in K^{*2}.$
\item Define $A := \OO_L $, the ring of integers  of $L$   if $a \notin K^{*2}.$ If $ a \in K^{*2}$, then  set $A:= \OO_K \times \OO_K$. 
\item Next, we define the quadratic \'etale $K_\q$-algebras: $L_\q=\frac{K_\q[X]}{(X^2-a)}$. 
Thus  $L_\q/K_\q$ is a quadratic field extension, if $ a \notin K_\q^{*2}$; otherwise, $L_\q \cong K_\q \times K_\q, \text{ if } a \in K_\q^{*2}$. 
\item Define $A_\q:=\begin{cases} \text{the ring of integers $\OO_{L_\q}$ of the field } L_\q, & \text{ if } a \notin K_\q^{*2},\\  \OO_{K_\q} \times \OO_{K_\q}, & \text{ if } a \in K_\q^{*2}. \end{cases}$ 
\item $N_{L/K}:L \to K$ denotes the field norm if $L$ is field and the multiplication of co-ordinates if $L \cong K \times K$. A similar notation is followed for $L_\q$ and $K_\q$.
\end{itemize} 

For $F \in \{K, K_\q\}$, let $\delta_{\psi_{a,b}, F}:\widehat{E}_{a,b}(F) \longrightarrow \widehat{E}_{a,b}(F)/\psi_{a,b}(E(F)) \lhook\joinrel\xrightarrow{\overline{\delta}_{\psi_{a,b}, F}} H^1(G_F, E[\psi_{a,b}])$ be the Kummer map. Then we have the following commutative diagram:
\begin{center}\begin{tikzcd}
 0 \arrow[r] & \widehat{E}_{a,b}(K)/\psi_{a,b}(E_{a,b}(K)) \arrow[r, "\overline{\delta}_{\psi_{a,b}, K}"] \arrow[d]
& H^1(G_K, E_{a,b}[\psi_{a,b}]) \arrow[d, "\underset{\q \in \Sigma_K}{\prod} {\rm res}_\q"] \arrow[r] & H^1(G_K,E_{a,b})[\psi_{a,b}] \arrow[d] \arrow[r] & 0\\
 0 \arrow[r] & \underset{\q \in \Sigma_K}{\prod} \widehat{E}_{a,b}(K_\q)/\psi_{a,b}(E_{a,b}(K_\q)) \arrow[r, "\underset{\q \in \Sigma_K}{\prod} \overline{\delta}_{\psi_{a,b},K_\q}"] & \underset{\q \in \Sigma_K}{\prod} H^1(G_{K_\q},E_{a,b}[\psi_{a,b}]) \arrow[r] & \underset{\q \in \Sigma_K}{\prod} H^1(G_{K_\q},E_{a,b})[\psi_{a,b}] \arrow[r] & 0.
\end{tikzcd}\end{center}

\begin{defn}\label{mainsel}
The $\psi_{a,b}$-Selmer group of $E_{a,b}$ over $K$, denoted ${\rm Sel}^{\psi_{a,b}}(E_{a,b}/K)$, is defined as
$${\rm Sel}^{\psi_{a,b}}(E_{a,b}/K) := \big\{\xi \in H^1(G_K,E_{a,b}[\psi_{a,b}]) \mid {\rm res}_\q(\xi) \in {\rm Im } \ \delta_{\psi_{a,b},K_\q} \text{ for every } \q \in \Sigma_K \big\}.$$
\end{defn}
\noindent Setting $\Sh(E_{a,b}/K):=\Ker \big( H^1(G_K,E_{a,b}) \to \underset{\q \in \Sigma_K}{\prod} H^1(G_{K_\q},E_{a,b}) \big)$ to be the Tate-Shafarevich group of $E_{a,b}$ over $K$, we get the fundamental exact sequence:
\begin{equation}\label{eq:defofsha}
0 \longrightarrow {\widehat{E}_{a,b}(K)}/{\psi_{a,b}(E_{a,b}(K))} \longrightarrow {\rm Sel}^{\psi_{a,b}}(E_{a,b}/K) \longrightarrow \Sh(E_{a,b}/K)[\psi_{a,b}] \longrightarrow 0
\end{equation}
For the isogeny given by the multiplication-by-$3$ map on $E_{a,b}$, the $3$-Selmer group ${\rm Sel}^3(E_{a,b}/K)$ is defined in a similar manner. 

The following result can be deduced from \cite[Prop.~24]{bes}.
\begin{proposition}\label{cohom}
There is an isomorphism of group schemes 
$E_{a,b}[\psi_{a,b}] \cong \Ker( \mathrm{Res}_K^{L_a} \mu_3 \to \mu_3)$
and an induced isomorphism
$$H^1(G_K, E_{a,b}[\psi_{a,b}]) \cong (L_a^*/L_a^{*3})_{N=1},$$
where $(L_a^*/L_a^{*3})_{N=1}$ denotes the kernel of the norm map $\overline{N}_{L_a/K}: L_a^*/L_a^{*3} \to K^*/K^{*3}$. 

The same statement holds if we replace $(E_{a,b},\psi_{a,b})$ by $(\widehat{E}_{a,b},\widehat{\psi}_{a,b})$. \qed
\end{proposition}
We denote the composite map $\widehat{E}_{a,b}(K) \longrightarrow \widehat{E}_{a,b}(K)/\psi_{a,b}(E_{a,b}(K))  \lhook\joinrel\xrightarrow{} H^1(G_K, E_{a,b}[\psi_{a,b}]) \isomto (L^*/L^{*3})_{N=1}$ again by $\delta_{\psi_{a,b},K}$, in view of Prop. \ref{cohom}. Similarly, we will view $\delta_{\psi_{a,b},K_\q}(\widehat{E}_{a,b}(K_\q))$ as a subset of $(L_\q^*/L_\q^{*3}\big)_{N=1}$. 
We have the commutative diagram,
\begin{equation}\label{eq:localglobal}
	\begin{tikzcd}
	  {{L^*}/{L^{*3}}} \arrow[r, "\overline{N}_{L/K}"] \arrow[d] & {{K^*}/{K^{*3}}} \arrow[d]\\
		{{L_\q^*}/{L_\q^{*3}}} \arrow[r, "\overline{N}_{L_\q/K_\q}"] & {{K_\q^*}/{K_\q^{*3}.}}
	\end{tikzcd}
 \end{equation}
  To ease the notation, via the canonical embedding $\iota_\q:K \to K_\q$, we identify $\iota_\q(x)$ with $x \in K$; similarly for $L$. Hence from \eqref{eq:localglobal}, if $\overline{x} \in \big(L^*/L^{*3}\big)_{N=1}$, then $\overline{x} \in \big(L_\q^*/L_\q^{*3}\big)_{N=1}$.

We now have the following alternative description of ${\rm Sel}^{\psi_{a,b}}(E_{a,b}/K)$:
\begin{equation}\label{eq:newseldef}
{\rm Sel}^{\psi_{a,b}}(E_{a,b}/K)=\big\{ \overline{x} \in (L^*/L^{*3})_{N=1} \mid \overline{x} \in \text{Im } \delta_{\psi_{a,b},K_\q} \text{ for all } \q \in \Sigma_K \big\}.
\end{equation}
Similarly, if $a \in K^{*2}$, then the definition in \eqref{eq:newseldef} can be written more explicitly:
\begin{equation}\label{eq:newseldefsq}
{\rm Sel}^{\psi_{a,b}}(E_{a,b}/K)=\big\{ (\overline{x}_1,\overline{x}_2) \in (K^*/K^{*3} \times K^*/K^{*3})_{N=1} \mid (\overline{x}_1,\overline{x}_2) \in \text{Im } \delta_{\psi_{a,b},K_\q} \text{ for all } \q \in \Sigma_K \big\}.
\end{equation}
These will be our working definitions of the Selmer groups for the rest of the article.

\subsection{Some algebraic number theory}\label{ant}
\noindent Here, we will recall some basic algebraic results from \cite{jms} which will be used later. The proofs can be found in \cite[\S2]{jms}.
\begin{defn}\label{defnofAq}
We define $(A^*/A^{*3})_{N=1}$ to be  the kernel of the norm map $\overline{N}_{L/K}: A^*/A^{*3} \to \OO_K^*/\OO_K^{*3}$ induced from the norm map $\overline{N}_{L/K}: L^*/L^{*3} \to K^*/K^{*3}$. The group 
$(A_\q^*/A_\q^{*3})_{N=1}$ is defined similarly by replacing $A$ and $\OO_K$ with $A_\q$ and $\OO_{K_\q}$, respectively, in this definition. 
\end{defn}

\begin{proposition}\label{N1subgrp}
For $\q \nmid 3$, we have
$\Big|\Big(\frac{A_\q^*}{A_\q^{*3}}\Big)_{N=1}\Big| =       \begin{cases}
1, & \text{ if } a \notin K_\q^{*2}, \\
3, & \text{ if } a \in K_\q^{*2}.   \end{cases}$\\
On the other hand, for $\p \mid 3$, we have 
$\Big|\Big(\frac{A_\p^*}{A_\p^{*3}}\Big)_{N=1}\Big| =       \begin{cases}
9, & \text{ if } a \notin K_\p^{*2}, \\
27, & \text{ if } a \in K_\p^{*2}. \end{cases}$ \qed
\end{proposition}

For a local field $F$ with the ring of integers $\OO_F$ and uniformizer $\pi$, we let $U^n_F:=1+\pi^n\OO_F$. Given an element $u \in \OO_F^*$, we let $\overline{u}$ denote its image in ${\OO_F^*}/{U^n_F} \cong (\OO_F/{\pi^n\OO_F})^*$.

\begin{defn}\label{defnofV3}
For a local field $F$, define
$$V_F = \big\{ u \in \OO_F^* \mid \overline{u}  \in  {\OO_F^*}/{U^3_F} \text{ satisfies } \overline{u} = \alpha^3 \text{ for some } \alpha \in {\OO_F^*}/{U^3_F} \big\}.$$

Next, we define a subset $V_3$ of $A_\p^*$:
$$V_3=\begin{cases} \big\{u \in V_{L_\mathfrak{p}} \mid N_{L_\mathfrak{p}/K_\p}(u) \in \OO_{K_\p}^{*3} \big\}, & \text{  if  } a \notin K_\p^{*2},\\
\big\{ (u_1,u_2) \in V_{K_\p} \times V_{K_\p}  \mid u_1u_2 \in \OO_{K_\p}^{*3} \big\}, & \text{  if  } a \in K_\p^{*2}. \end{cases}$$
\end{defn}
Note that $A_\p^{*3} \subseteq V_3$ and $V_3/A_\p^{*3}  \subseteq (A_\p^*/A_\p^{*3})_{N=1}$, whether or not $a$ is a square in $K_\p^*$. 

\begin{lemma}\label{V3U3U4}
Let $a \in K_\p^{*2}$, then $\Big(\frac{A_\p^*}{A_\p^{*3}}\Big)_{N=1} \cong \frac{U^1_{K_\p}}{U^4_{K_\p}}$ and $\frac{V_3}{A_\p^{*3}} \cong \frac{U^3_{K_\p}}{U^4_{K_\p}}$. \qed
\end{lemma}

\begin{proposition}\label{sizeofV3}
We have $\Big|\frac{V_3}{A_\p^{*3}}\Big| = 1$ if $a \notin K_\p^{*2}$ and $\Big|\frac{V_3}{A_\p^{*3}}\Big| =3$ if $a \in K_\p^{*2}$. \qed
\end{proposition}

\noindent {\bf Notation}: For a non-archimedean local field $F$ with discrete valuation $v$ and $\overline{x} \in {F}^*/{{F}^{*3}}$,  define $v(\overline{x}) \in \Z/{3\Z}$ by $v(\overline{x}):= v(x) \pmod 3$, for  any lift $x \in F^*$. 
Similarly, for an element $(\overline{y},\overline{z})$ of $F^*/{F^{*3}} \times F^*/{F^{*3}}$, we define $v(\overline{y}, \overline{z}):=\big(v(y) \pmod 3, v(z) \pmod 3 \big) \in \Z/{3\Z} \times \Z/3\Z$. If both $v(y), v(z) \equiv 0 \pmod 3$,  then by abuse of notation, we write $v(\overline{y},\overline{z}) \equiv 0 \pmod{3}$.

For a prime $\q \in \Sigma_K$, we denote the induced valuation on $L_\q$ simply by $\upsilon_\q$.

\begin{lemma}\label{norm1mult3}
For an element $\overline{x}$ of $\Big(\frac{L_\q^*}{L_\q^{*3}}\Big)_{N=1}$, one has $\overline{x} \in \Big(\frac{A_\q^*}{A_\q^{*3}}\Big)_{N=1}$ if and only if $\upsilon_\q(\overline{x}) \equiv 0 \pmod 3$. \qed
\end{lemma}

\begin{defn}\label{defofM1M2}
Let $S$ be a finite set of finite primes of $K$ and recall that $\OO_S$ denotes the ring of $S$-integers in $K$.
If $a \notin K^{*2}$, then $L=K(\sqrt{a})$. Let $S(L)=\{ \mathfrak{Q} \in \Sigma_L \mid \mathfrak{Q} \cap K \in S \}$. Further, $\OO_{S(L)}$ will denote the ring of $S(L)$-integers in $L$. Let $x \in L^*$ be any lift of $\overline{x} \in L^*/L^{*3}$. We define the following $\F_3$-vector spaces:
$$M(S,a)=\big\{ \overline{x} \in L^*/L^{*3} : L(\sqrt[3]{x})/L \text{ is unramified and } x \in L_\q^{*3} \text{ for all } \q  \in  S \big\},$$
$$N(S,a)= \big\{ \overline{x} \in L^*/L^{*3} : (x)=I^3 \text{ for some fractional ideal } I \text{ of } \OO_{S(L)}  \big\}.$$
\end{defn}

Note that if $\overline{x} \in M(S,a)$, then any lift $x \in L^*$ of $\overline{x}$ satisfies $(x)=I^3$ for some fractional ideal $I$ of $\OO_{S(L)}$ (\cite[Ch.~9]{was}). This implies that $N_{L/K}(x) \in K^{*3}$. So, $M(S,a) \subset N(S,a) \subset (L^*/L^{*3})_{N=1}$.

\begin{defn}
Let $S$ be a finite set of finite primes of $K$. When $a \in K^{*2}$, we define 
$$N'(S,a):=\big\{ (\overline{x}_1,\overline{x}_2) \in \big(K^*/K^{*3} \times K^*/K^{*3}\big)_{N=1} \mid \upsilon_\q(\overline{x}_1,\overline{x}_2) \equiv 0 \pmod 3 \ \text{ for all } \ \q \notin S \big\}.$$
\end{defn}
\begin{defn}\label{defofh3}
We assume that  $a \notin K^{*2}$ and denote by $Cl_{S(L)}(L)$ the $S(L)$-ideal class group of $L$ i.e. the class group of $\OO_{S(L)}$. Further, let $h^3_{S(L)}$ denote the $3$-rank of $Cl_{S(L)}(L)$ i.e. $h^3_{S(L)}:= \dim_{\F_3} Cl_{S(L)}(L) \otimes_\Z \F_3 = \dim_{\F_3}Cl_{S(L)}(L)[3]$. In particular, if $S(L)$ is empty, then $Cl_\emptyset(L)=Cl(L)$ is the ideal class group of $L$ and $h^3_L=\dim_{\F_3} Cl(L) \otimes_\Z \F_3$.
\end{defn}

\begin{proposition}\label{boundsgen}
If $a \notin K^{*2}$, we have $\dim_{\F_3} M(S,a) = h^3_{S(L)}$ and $\dim_{\F_3} N(S,a)= h^3_{S(L)}+ |S(L)|+2$.

On the other hand, one has that $\dim_{\F_3} N'(S,a) = |S|+1$ when $a \in K^{*2}$. \qed
\end{proposition}

We have an equivalently description of the $\F_3$-modules $M(S,a)$, $N(S,a)$ and $N'(S,a)$, which are compatible with the definition of Selmer groups given in \eqref{eq:newseldef} and \eqref{eq:newseldefsq}.
\begin{proposition}\label{newM1M2sel}
When $a \notin K^{*2}$, we have
$$M(S,a)=\Big\{ \overline{x} \in \Big(\frac{L^*}{L^{*3}}\Big)_{N=1} \Big\vert \ \overline{x} \in \Big(\frac{A_\q^*}{A_\q^{*3}}\Big)_{N=1} \text{ if } \q \notin S \cup \{\p\}, \text{ } \overline{x}=\overline{1} \text{ if } \q \in S \text{ and } \overline{x} \in \frac{V_3}{A_\p^{*3}} \text{ if } \p \notin S \Big\}$$
and $N(S,a)=\Big\{ \overline{x} \in \Big(\frac{L^*}{L^{*3}}\Big)_{N=1} \Big\vert \ \overline{x} \in \Big(\frac{A_\q^*}{A_\q^{*3}}\Big)_{N=1} \text{ for all } \q \notin S \Big\}.$

\noindent Whereas for $a \in K^{*2}$, $N'(S,a)=\Big\{ (\overline{x}_1,\overline{x}_2) \in \Big(\frac{K^*}{K^{*3}} \times \frac{K^*}{K^{*3}}\Big)_{N=1} \Big\vert \ (\overline{x}_1,\overline{x}_2) \in \Big(\frac{A_\q^*}{A_\q^{*3}}\Big)_{N=1} \text{ for all } \q \notin S \Big\}$. \qed
\end{proposition}

\section{Local theory for the curves $E_{a,b}$}\label{type2curves}
Let $F$ be a local field of characteristic $0$ and $\omega$ be a uniformizer of  $\OO_F$. Suppose $\varphi: E \to \widehat{E}$ be an isogeny of elliptic curves and $c_\omega(E)$ and $c_\omega(\widehat{E})$ denote the Tamagawa numbers of $E$ and $\widehat{E}$ at  $\omega$. By a theorem of Schaefer \cite[Lemma~3.8]{sch}, we know that
\begin{equation}\label{eq:schaeferfor}
|\widehat{E}(F)/\varphi(E(F))|  = \frac{||\varphi'(0)||_\omega^{-1} \times |E(F)[\varphi]| \times c_\omega(\widehat{E}) }{c_\omega(E)},
\end{equation}
where $||.||_\omega$ is the norm in $F$ and $\varphi'(0)$ is the leading coefficient of the power series representation of $\varphi$ on the formal group of $\widehat{E}$ \cite[\S4]{sil1}.
Using \eqref{eq:schaeferfor}, we will compute $|\widehat{E}(K_\q)/\varphi(E(K_\q))|$ for all $\q$.\\

Consider the curves $E_{a,b}: y^2=x^3+a(x-b)^2$, with $ab(4a+27b) \neq 0$, $a, b \in \Z$ such that $gcd(a,b)$ is square-free in $K$ and the rational $3$-isogeny  $\psi_{a,b}: E_{a,b} \to \widehat{E}_{a,b}$,
where $\widehat{E}_{a,b}: y^2=x^3-27a(x-d)^2$ with $d=4a+27b$. In this subsection, we give an explicit description of the Kummer maps $\delta_{\psi_{a,b},K_\q}$ for all $\q \in \Sigma_K$.
But first, we state a general result (cf. \cite[Lemma~3.1]{jms} for the proof).

\begin{lemma}\label{goodreduction}
Let $\varphi: E \to \widehat{E}$ be a $3$-isogeny defined over $K$. Recall from \S \ref{iso} that $\delta_{\varphi,K_\q}(E(K_\q)) \subset (L_\q^*/L_\q^{*3})_{N=1}$.
If $E$ (and hence $\widehat{E}$) has good reduction at a prime $\q \nmid 3$ in $K$, then for any $P \in \widehat{E}(K_\q)$, \ $v_\q(\delta_{\varphi,K_\q}(P)) \equiv 0 \pmod 3$. In particular, $\delta_{\varphi,K_\q}(\widehat{E}(K_\q)) \subset \big(A_\q^*/A_\q^{*3}\big)_{N=1}$. \qed
\end{lemma}

For $P \in \widehat{E}(K_\q)$, $\delta_{{\varphi},K_\q}(P) = \overline{t} \in L_\q^*/L_\q^{*3}$, if $L_\q$ is a field and $\delta_{{\varphi},K_\q}(P)=(\overline{t_1},\overline{t_2})$, if $L_\q \cong K_\q \times K_\q$. To ease the notation, we write $\delta_{{\varphi},K_\q}(P)=t$ or $(t_1,t_2)$, respectively.

\begin{proposition}\cite[\S14, \S15]{cass}\label{type2}
Let  $\psi_{a,b}: E_{a,b}  \to \widehat{E}_{a,b}$
be the $K$-isogeny defined in \eqref{eq:defofpsiab}. 
Recall that $d=4a+27b$.
We have the Kummer map $\delta_{{\psi}_{a,b},K_\q}: \widehat{E}_{a,b}(K_\q) \to  \big(L_\q^*/L_\q^{*3}\big)_{N=1}$ with the kernel  $\Ker(\delta_{{\psi}_{a,b},K_\q})={\psi}_{a,b}({E}_{a,b}(K_\q))$.\\
\noindent If $L_\q$ is a field, then $\delta_{{\psi}_{a,b},K_\q}$  is given by $\delta_{{\psi}_{a,b},K_\q}(P) = \begin{cases}
{1},                            & \text{ if } P = O,\\
{y(P) - \p^3\sqrt{a}\big(x(P)- d\big)}, & \text{ otherwise.} \end{cases}$\\
\noindent If $L_\q \cong K_\q \times K_\q$, then $\delta_{{\psi}_{a,b},K_\q}$ is given by  
$$\delta_{{\psi}_{a,b},K_\q}(P) = \begin{cases}
({1},{1}),                                                        & \text{ if } P = O,\\
\big(2{\p^3\sqrt{a}d}, \ \frac{1}{2 {\p^3\sqrt{a} d}}\big),                 & \text{ if } P = (0, \p^3\sqrt{a}d),\\
\big(-\frac{1}{2 {\p^3\sqrt{a}d}}, \ -2{\p^3\sqrt{a}d}\big),                & \text{ if } P = (0, -\p^3\sqrt{a}d),\\
\big({y(P) - \p^3\sqrt{a}(x(P)- d)}, \ {y(P) + \p^3\sqrt{a}(x(P)- d)}\big), & \text{ if } P \notin \widehat{E}_{a,b}(K_\q)[\widehat{\psi}_{a,b}]. \qed
\end{cases}$$
\end{proposition}

\begin{lemma}{\cite[\S~7]{sil1}}\label{redntype}
Let $\ell \in \Z$ be a prime. If $F$ is a local field, $\omega$ is a prime in $F$ and $\varphi: E \to \widehat{E}$ is an $\ell$-isogeny over $F$, then
\begin{enumerate}
\item $E$ has good reduction at $\omega$ if and only if $\widehat{E}$ has good reduction at $\omega$.
\item $E$ has split multiplicative reduction at $\omega$ if and only if $\widehat{E}$ has split multiplicative reduction at $\omega$. In this situation, one of $E$ and $\widehat{E}$ is of Kodaira type $I_n$, while the other is of the type $I_{\ell n}$ with $n=-\upsilon_\omega(j(E))$.\qed
\end{enumerate} 
\end{lemma}

\subsection{\bf{Local theory at primes $\q \nmid 3$}}
Let $\q$ be a prime of $K=\Q({\zeta_3})$ not dividing $3$.
\begin{proposition}\label{type2nsq}
If $a \notin K_\q^{*2}$, then we have 
 $\frac{\widehat{E}_{a,b}(K_\q)}{\psi_{a,b}(E_{a,b}(K_\q))} =\Big(\frac{A_\q^*}{A_\q^{*3}}\Big)_{N=1} = \{1\}.$
\end{proposition}

\begin{proof}
Since $a \notin K_\q^{*2}$, we have $E_{a,b}(K_\q)[\psi_{a,b}]=\{O\}$. 
Using \eqref{eq:schaeferfor}, we get
$\Big|\frac{\widehat{E}_{a,b}(K_\q)}{\psi_{a,b} (E_{a,b}(K_\q))}\Big| = \frac{c_\q(\widehat{E}_{a,b})}{c_\q(E_{a,b})}$ $ \mbox{ and } \Big|\frac{E_{a,b}(K_\q)}{\widehat{\psi}_{a,b} (\widehat{E}_{a,b}(K_\q))}\Big| = \frac{c_\q(E_{a,b})}{c_\q(\widehat{E}_{a,b})}$.
Thus, $\frac{\widehat{E}_{a,b}(K_\q)}{\psi_{a,b}(E_{a,b}(K_\q))}$ is trivial, also $\big({A_\q^*}/{A_\q^{*3}}\big)_{N=1} = \{1\}$ by Prop. \ref{N1subgrp}.
\end{proof}

\begin{proposition}\label{type2redn}
 Suppose $a \in K_\q^{*2}$ and $E_{a,b}$ (or equivalently $\widehat{E}_{a,b}$) has good reduction at $\q$. Then we have that 
$\big|\frac{\widehat{E}_{a,b}(K_\q)}{\psi_{a,b}(E_{a,b}(K_\q))}\big| = 3$ and
hence $\frac{\widehat{E}_{a,b}(K_\q)}{\psi_{a,b}(E_{a,b}(K_\q))} \cong \Big(\frac{A_\q^*}{A_\q^{*3}}\Big)_{N=1}.$
\end{proposition}

\begin{proof}
As $a \in K_\q^{*2}$, $\big|E_{a,b}[\psi_{a,b}](K_\q)\big|=3$ and again $||\psi_{a,b}'(0)||_\q^{-1}=1$ by \cite{sch}. 
Since $\q$ is a good prime for $E_{a,b}$ (and $\widehat{E}_{a,b}$), we have $\big|\frac{\widehat{E}_{a,b}(K_\q)}{\psi_{a,b}(E_{a,b}(K_\q))}\big| =3$ by \eqref{eq:schaeferfor}. The result follows from Prop. \ref{N1subgrp} and Lemma \ref{goodreduction}.
\end{proof}

Next, let $\q \nmid 3$ be a prime of $K$ such that  $a \in K_\q^{*2}$ and $\q \mid \Delta_{E_{a,b}}=-2^4a^2b^3d$. The local data obtained for each such $\q$ 
using Tate's algorithm \cite[Ch.~4, Algo.~9.4]{sil2}, the Schaefer's formula \eqref{eq:schaeferfor} and Prop. \ref{N1subgrp} is consolidated in Table \ref{tab:type2atq} below.

\begin{table}[h]
\centering
	\begin{tabular}{ |c||c|c| }
		\hline
		 $a \in K_\q^{*2}$ and & \multirow{2}{*}{Tamagawa Numbers} &  \\
		 $\q$ divides $\Delta_{E_{a,b}}$ & &  \\
		\hline\hline
		$\q \mid a$ and & \multirow{2}{*}{$c_\q(E_{a,b})=c_\q(\widehat{E}_{a,b})=1$} & \multirow{2}{*}{$\frac{\widehat{E}_{a,b}(K_\q)}{\psi_{a,b}(E_{a,b}(K_\q))} = \Big(\frac{A_\q^*}{A_\q^{*3}}\Big)_{N=1}$}\\
		$\upsilon_\q(4ab^2) \equiv 0 \pmod 6$  & &   \\
		\hline
		$\q \mid a$ and & \multirow{2}{*}{$c_\q(E_{a,b})=c_\q(\widehat{E}_{a,b})=3$} & \multirow{2}{*}{$\frac{\widehat{E}_{a,b}(K_\q)}{\psi_{a,b}(E_{a,b}(K_\q))} \cap \Big(\frac{A_\q^*}{A_\q^{*3}}\Big)_{N=1}=\{1\}$}\\
		$\upsilon_\q(4ab^2) \equiv 2 \pmod 6$  & &  \\
		\hline
		$\q \mid a$ and & \multirow{2}{*}{$c_\q(E_{a,b})=c_\q(\widehat{E}_{a,b})=3$} & \multirow{2}{*}{$\frac{\widehat{E}_{a,b}(K_\q)}{\psi_{a,b}(E_{a,b}(K_\q))} \cap \Big(\frac{A_\q^*}{A_\q^{*3}}\Big)_{N=1}=\{1\}$}\\
		$\upsilon_\q(4ab^2) \equiv 4 \pmod 6$  & & \\
		\hline
		\multirow{2}{*}{$\q \mid b$ and $\q \nmid 2a$} & {$c_\q(E_{a,b})=3\upsilon_\q(b)$} & \multirow{2}{*}{$\{1\}=\frac{\widehat{E}_{a,b}(K_\q)}{\psi_{a,b}(E_{a,b}(K_\q))} \subsetneq \Big(\frac{A_\q^*}{A_\q^{*3}}\Big)_{N=1}$}\\
		& {$c_\q(\widehat{E}_{a,b})=\upsilon_\q(b)$} & \\
		\hline
		\multirow{2}{*}{$\q \mid d$ and $\q \nmid 2a$} & {$c_\q(E_{a,b})=\upsilon_\q(d)$} & \multirow{2}{*}{$\Big(\frac{A_\q^*}{A_\q^{*3}}\Big)_{N=1} \subsetneq \frac{\widehat{E}_{a,b}(K_\q)}{\psi_{a,b}(E_{a,b}(K_\q))}$}\\
		& {$c_\q(\widehat{E}_{a,b})=3\upsilon_\q(d)$} & \\
		\hline
		$\q=2$ & \multirow{2}{*}{$c_\q(E_{a,b})=c_\q(\widehat{E}_{a,b})=3$} & \multirow{2}{*}{$\frac{\widehat{E}_{a,b}(K_\q)}{\psi_{a,b}(E_{a,b}(K_\q))} \cap \Big(\frac{A_\q^*}{A_\q^{*3}}\Big)_{N=1}=\{1\}$}\\
		$2 \nmid a$ and $2 \nmid b$  & &  \\
		\hline
		$\q=2$ & \multirow{2}{*}{$c_\q(E_{a,b})=c_\q(\widehat{E}_{a,b})=3$} & \multirow{2}{*}{$\frac{\widehat{E}_{a,b}(K_\q)}{\psi_{a,b}(E_{a,b}(K_\q))} \cap \Big(\frac{A_\q^*}{A_\q^{*3}}\Big)_{N=1}=\{1\}$}\\
		$2 \nmid a$ and $2 || b$  & &  \\
		\hline
		$\q=2$, & {$c_\q(E_{a,b})=\upsilon_\q(d)-2$} & \multirow{2}{*}{$\Big(\frac{A_\q^*}{A_\q^{*3}}\Big)_{N=1} \subsetneq \frac{\widehat{E}_{a,b}(K_\q)}{\psi_{a,b}(E_{a,b}(K_\q))}$}\\
		$2 \nmid a$ and $4 || b$ & {$c_\q(\widehat{E}_{a,b})=3(\upsilon_\q(d)-2)$}  & \\
		\hline
		$\q=2,$  & {$c_\q(E_{a,b})=3(\upsilon_\q(b)-2)$} & \multirow{2}{*}{$\{1\}=\frac{\widehat{E}_{a,b}(K_\q)}{\psi_{a,b}(E_{a,b}(K_\q))} \subsetneq \Big(\frac{A_\q^*}{A_\q^{*3}}\Big)_{N=1}$ }\\
		$2 \nmid a$ and $8 \mid b$ & {$c_\q(\widehat{E}_{a,b})=\upsilon_\q(b)-2$} & \\
		\hline
	\end{tabular}
\caption{Local data at a prime $\q\nmid 3$}
\label{tab:type2atq}
\end{table}

\subsection{\bf{Local theory at the prime $\p \mid  3$} }
Recall that $\p$ is the unique prime dividing $3$ in $K$.
\begin{proposition}\label{ais-1modulo3}
For $a \notin K_\p^{*2}$, one has $\frac{\widehat{E}_{a,b}(K_\p)}{\psi_{a,b}(E_{a,b}(K_\p))} \subset \Big(\frac{A_\p^*}{A_\p^{*3}}\Big)_{N=1}$.
\end{proposition}

\begin{proof}
The proof is similar to \cite[Lemma~3.9]{jms}.
\end{proof}

\noindent We divide the study of image of local Kummer maps $\delta_{\psi_{a,b},K_\p}$ at $\p$ into three subcases $3 \nmid ab$, $3 \mid b$ and $3 \mid a$.
\subsubsection{\underline {{\bf Case-1: $3 \nmid ab$ }} }\label{case1type2}
In this case, we find a $\p$-minimal model for $\widehat{E}_{a,b}$ 
and use it to compute the image of the Kummer map $\delta_{\psi_{a,b},K_\p}$. 

Note that $E_{a,b}:y^2=x^3+a(x-b)^2$ is $\p$-minimal and has good reduction at $\p$.
By Tate's algorithm \cite{sil2}, the $\p$-minimal model for $\widehat{E}_{a,b}$ is given by $E':=E'_{a,b}$
\begin{equation}\label{eq:defnofE2}
E' : Y^2= X^3+ a_2 X^2 + b_2 X + c_2,
\end{equation}
where $a_2 = -a(3+ 5\zeta_3)$, $b_2 = a \big( a(-5+ 2 \zeta_3) +18b\zeta_3^2 \big)$ and $c_2 = a\big( a^2(2+3\zeta_3)-2ab(1-9\zeta_3)-27b^2 \big).$

The following map $\theta: E' \to E_{-27a,d}$ is then a $G_K$-isomorphism
\begin{equation}\label{eq:defoftheta}
\theta(x(P), y(P)) = (x'(P),y'(P)) := ( \p^2(\p^2x(P) -a - a \p), \p^6 y(P)).
\end{equation}

The isogeny $\Psi_{a,b}:=\theta^{-1}\circ\psi_{a,b} : E_{a,b} \to E'$ 
is given by

\begin{equation}\label{eq:newpsi}
 \Psi_{a,b}(x,y) = (X,Y):= \Big( \frac{ \zeta_3( x^3 - a(3\p+5/3)x^2 - 4abx + 4ab^2 ) }{x^2}, \frac{-y ( x^3 + 4abx - 8ab^2 )}{x^3} \Big).
\end{equation}

Then, using formal groups of elliptic curves as in \cite[pg.~92]{sch}, we compute
$$Z:=-\frac{X}{Y} = \zeta_3 \cdot \frac{x}{y} \Big( 1 - \frac{ a(3\p+5/3)x^2 + 8abx - 12ab^2}{x^3 + 4abx - 8ab^2} \Big) = -{\zeta_3} \cdot z + O(z^2).$$
Hence, $\Psi'_{a,b}(0)=-{\zeta_3}$ and $||\Psi'_{a,b}(0)||_\p^{-1}=1$.

\begin{proposition}\label{3nmidab1}
Assume that $3 \nmid ab$. We have
\begin{enumerate}
    \item If $a \notin K_\p^{*2}$, then $\Big|\frac{\widehat{E}_{a,b}(K_\p)}{\psi_{a,b}(E_{a,b}(K_\p))}\Big| =1$ and $\frac{\widehat{E}_{a,b}(K_\p)}{\psi_{a,b}(E_{a,b}(K_\p))} \subsetneq \Big(\frac{A_\p^*}{A_\p^{*3}}\Big)_{N=1}.$
		
    \item If $a \in K_\p^{*2}$, then $\Big|\frac{\widehat{E}_{a,b}(K_\p)}{\psi_{a,b}(E_{a,b}(K_\p))}\Big| =3$ and $\frac{\widehat{E}_{a,b}(K_\p)}{\psi_{a,b}(E_{a,b}(K_\p))} \subsetneq \Big(\frac{A_\p^*}{A_\p^{*3}}\Big)_{N=1}$.
\end{enumerate}
\end{proposition}

\begin{proof}
We know that $c_\p(E_{a,b})=1$, as $E_{a,b}$ has good reduction at $\p$; thus $c_\p(E')=1$. Also, we have $||\Psi_{a,b}'(0)||_{\p}^{-1}=1$.
Using the Schaefer's formula \eqref{eq:schaeferfor}, we obtain $\big| \frac{{E}'(K_\p)}{\Psi_{a,b}(E_{a,b}(K_\p))} \big| = \begin{cases}
    1, & \text{ if } a \notin K_\p^{*2},\\
    3, & \text{ if } a \in K_\p^{*2}.\end{cases}$
Now, $\big| \frac{\widehat{E}_{a,b}(K_\p)}{\psi_{a,b}(E_{a,b}(K_\p))} \big|=\big| \frac{\theta({E}'(K_\p))}{\theta \circ \Psi_{a,b}(E_{a,b}(K_\p))} \big|=\big| \frac{{E}'(K_\p)}{\Psi_{a,b}(E_{a,b}(K_\p))} \big|$ 
and so, $\big| \frac{\widehat{E}_{a,b}(K_\p)}{\psi_{a,b}(E_{a,b}(K_\p))} \big| = \begin{cases}
    1, & \text{ if } a \notin K_\p^{*2},\\
    3, & \text{ if } a \in K_\p^{*2}.\end{cases}$

When $a \notin K_\p^{*2}$, then $\Big|\Big(\frac{A_\p^*}{A_\p^{*3}}\Big)_{N=1}\Big|=9$  by Prop. \ref{N1subgrp} 
and the strict inequality in (1) follows. When $a \in K_\p^{*2}$, then $\Big|\Big(\frac{A_\p^*}{A_\p^{*3}}\Big)_{N=1}\Big|=27$  by Prop. \ref{N1subgrp}. Thus, to complete the proof, it suffices to show, by Lemma \ref{norm1mult3}, that $\upsilon_\p(\delta_{\psi_{a,b},K_\p}(P')) \equiv 0 \pmod 3$ for all $P' \in \widehat{E}_{a,b}(K_\p)$. We now prove this by showing that for any $P \in E'(K_\p)$ we have $\upsilon_\p(\delta_{\psi_{a,b},K_\p}(\theta(P))) \equiv 0 \pmod 3$.

Consider a point $P \in {E}'(K_\p)[\widehat{\Psi}_{a,b}]$ and recall that $3 \nmid ab$ and $\p \mid 3$. Hence, $\upsilon_\p(a)= \upsilon_\p(d)=0$ and $\upsilon_\p(2\p^3 d\sqrt{a})=3$. Now, note here that $\theta(P)=(0, \pm \p^3 d\sqrt{a}) \in \widehat{E}_{a,b}(K_\p)[\widehat{\psi}_{a,b}]$ and we have $\upsilon_\p(\delta_{\psi_{a,b},K_\p}(\theta(P))) \equiv 0 \pmod 3$.
         
Next, take $P \notin E'(K_\p)[\widehat{\Psi}_{a,b}]$ and $\theta(P)=(x'(P),y'(P))$. If $\upsilon_\p(y'(P)) \neq \upsilon_\p(\p^3 \sqrt{a}(x'(P)-d))$, then as in the proof of \cite[Prop.~3.7(1)]{jms}, we get $\upsilon_\p(\delta_{\psi_{a,b},K_\p}((x'(P), y'(P)))) \equiv 0 \pmod 3$. So, we only consider the case $\upsilon_\p(y'(P)) = \upsilon_\p(\p^3 \sqrt{a}(x'(P)-d))=n$, say.

Further, if $\upsilon_\p(x'(P))=0$,
then $\upsilon_\p(y'(P)) =0$ but $\upsilon_\p(\p^3 \sqrt{a} (x'(P) -d)) \ge 3$, a contradiction.
Now if $\upsilon_\p(x'(P)) <0$, then  $x'(P)=\p^2 (\p^2x(P) - a - \p a)$ implies that $\upsilon_\p(x(P)) < 0$. This gives   $\upsilon_\p(x(P))= -2r$ for some $r \ge 3$.  
Thus, $\upsilon_\p(\p^3 \sqrt{a}(x'(P)-d))= -2r+7 \neq -3r+6 = \upsilon_\p(y'(P))$. This is a contradiction.
				
Finally, the case $x'(P) \neq 0$ and $\upsilon_\p(x'(P))>0$ and $\upsilon_\p(y'(P)) = \upsilon_\p(\p^3 \sqrt{a}(x'(P)-d))$ can be proved in a similar manner as in \cite[Lemma~3.9]{jms}. 
\end{proof}

\begin{proposition}\label{3nmidab2}
For $3 \nmid ab$ and $a \in K_\p^{*2}$, we have $\frac{\widehat{E}_{a,b}(K_\p)}{{\psi_{a,b}}(E_{a,b}(K_\p))} = \frac{V_3}{A_\p^{*3}}$. 
\end{proposition}
\proof
Note that $|\text{Im } (\delta_{\psi_{a,b},K_\p} \circ \theta)| =3$ by Prop. \ref{3nmidab1}. Let $P=(x(P),y(P)) \in E'(K_\p)$, $\theta(P)=(x'(P),y'(P))$ and $(C,D):=\delta_{\psi_{a,b},K_\p} \circ \theta(P)=\delta_{\psi_{a,b},K_\p}(x'(P),y'(P))$. Then $(C,D) \in (A_\p^*/A_\p^{*3})_{N=1}$ by Prop. \ref{3nmidab1}. Recall from Definition \ref{defnofV3} that $V_3 = \{ (u_1,u_2) \in V_{K_\p} \times V_{K_\p} \mid u_1u_2 \in \OO_{K_\p}^{*3} \}$ and $|V_3/A_\p^{*3}|=3$ from Prop. \ref{sizeofV3}.
Thus, it reduces to show that $\text{Im }(\delta_{\psi_{a,b},K_\p} \circ \theta) \subset V_3/A_\p^{*3}$.

By Lemma \ref{norm1mult3}, $(\upsilon_\p(C), \upsilon_\p(D)) \equiv 0 \pmod 3$. Also, as elements of $K^*/K^{*3} \times K^*/K^{*3}$,
$$(\overline{C}, \overline{D})=(\overline{\p^{-\upsilon_\p(C)}u_1^3C}, \ \overline{\p^{-\upsilon_\p(D)}u_2^3D}) \ \text{ for any } \ u_1, u_2  \in \OO_{K_\p}^*.$$
This gives $\p^{-\upsilon_\p(C)}u_1^3C \cdot \p^{-\upsilon_\p(D)}u_2^3D \in \OO_{K_\p}^{*3}.$ Hence, it is enough to show that there exists some $u \in \OO_{K_\p}^*$ such that either $\overline{\p^{- \upsilon_\p(C)}u^3 C} \in  (\OO_{K_\p}^*/U_{K_\p}^3)^{3}$ or $\overline{\p^{-\upsilon_\p(D)}u^3 D} \in  (\OO_{K_\p}^*/U_{K_\p}^3)^{3}$.

\noindent For $(x'(P), y'(P)) \notin \widehat{E}_{a,b}(K_\p)[\widehat{\psi}_{a,b}]$, recall that $C= y'(P) - \p^3 \sqrt{a} (x'(P)-d)$ and $D= y'(P) + \p^3 \sqrt{a} (x'(P)-d)$. 

If $\upsilon_\p(x(P)) \ge 0$, then  $\upsilon_\p(y(P)) \ge 0$. This implies, $\upsilon_\p(y'(P)) \ge 6$ and $\upsilon_\p(x'(P))=2$, which shows that $\upsilon_\p(C)=\upsilon_\p(D)=3$. Then 
$$\p^{-3}C = \p^3 y(P) - \sqrt{a}(\p^4 x(P) - \p^2a - \p^3a -d) \equiv a\sqrt{a} \p^2 + \sqrt{a} d \equiv a\sqrt{a} (\p^2 +4)  \equiv a\sqrt{a}\pmod{\p^3 },$$
where the last two congruences use the facts that $d \equiv 4a \pmod{\p^3}$ and $\p^2 + 4 \equiv 1 \pmod{\p^3}$. Hence, we obtain that 
$\overline{\p^{-\upsilon_\p(C)}(\sqrt{a})^{-3}C}=\overline{1} \in \big({\OO_{K_\p}^*}/{U_{K_\p}^3}\big)^{3}$.

For the rest of the proof, we deal with the case $\upsilon_\p(x(P)) < 0$. Thus, $\upsilon_\p(x(P))=-2r$ and $\upsilon_\p(y(P))=-3r$ for some $r \in \N$. By the theory of formal groups, there exists $0 \neq z = \p^r w$ with $w \in \OO_{K_\p}^*$ such that $P(z) = P = (x(P),y(P)):=(x(P(z)), y(P(z))) \in E'(K_\p)$. We divide this in two cases : $r \ge 2$ and $r=1$.

\begin{enumerate}
    \item Assume that $\upsilon_\p(x(P(z))) = -2r$, where $0 \neq z=\p^r w$ with $r \ge 2$.\\
Therefore, one has $\upsilon_\p(x'(P(z))) = -2r+4$ and $\upsilon_\p(y'(P(z))) = -3r+6$, which in turn implies that $\upsilon_\p(C)= \upsilon_\p(D) =-3r+6 \le 0$. Now,
\begin{align*}
\p^{3r-6}w^3C & = z^3 y(P(z)) - w\sqrt{a} \p^{r+1} z^2 x(P(z)) + w^3 a\sqrt{a} \p^{3r-1}+ w^3 a\sqrt{a} \p^{3r}+ d w^3 \sqrt{a} \p^{3r-3} \\
& \equiv z^3y(P(z)) \pmod{\p^3} \equiv -1 \pmod{\p^3}.
\end{align*}
Thus we have $\overline{\p^{-\upsilon_\p(C)}w^{-3}C}=-\overline{1} \in \big({\OO_{K_\p}^*}/{U_{K_\p}^3}\big)^{3}$ in this case.

\item  Now assume that $r=1$ i.e. $\upsilon_\p(x(P(z))) =-2$.
\begin{enumerate}
    \item First consider the case $x(P(z)) \neq \frac{1}{\p^2}(a+ \p a)$.\\
In this case, $\upsilon_\p(y'(P(z)))= 3$ and $\upsilon_\p(\p^3 \sqrt{a} (x'(P(z))-d)) =3$. Thus, at least one of $\upsilon_\p(C)$ and $\upsilon_\p(D)$ is $3$. Further, $w\sqrt{a} \equiv \pm 1 \pmod \p$ and we only prove the case $w\sqrt{a} \equiv -1 \pmod \p$ here. The case $w\sqrt{a} \equiv 1 \pmod \p$ can be proved in a similar manner.

Let  $z=\p w$ and $w\sqrt{a} \equiv -1 \pmod{\p}$. Then
\begin{align*}
\p^{-3}w^3C & \equiv z^3y(P(z)) -  w\sqrt{a} \p^2 z^2x(P(z)) +  w^3 a\sqrt{a} \p^2 +  dw^3 \sqrt{a} \pmod{\p^3}\\
& \equiv -1+ a_2\p^2 w^2 -  w \sqrt{a} \p^2 +  w^3 a\sqrt{a} \p^2 + 4  w^3 a\sqrt{a} \pmod{\p^3} \\
& \equiv -1+ 4  w^3 a\sqrt{a} + \p^2( a_2w^2+ w^3 a\sqrt{a} -  w\sqrt{a}) \pmod{\p^3}\\
& \equiv -5 + a_2 \p^2 \pmod{\p^3} \equiv -5 + a \p^2 \pmod{\p^3} \equiv -5+ \p^2 \pmod{\p^3} \equiv 1 \pmod{\p^3}. 
\end{align*}
Thus, we obtain $\overline{\p^{-\upsilon_\pi(C)}w^{3}C}=\overline{1} \in \big({\OO_{K_\p}^*}/{U_{K_\p}^3}\big)^{3}$ in this case.

\item Finally, consider the case $x(P(z)) = \frac{1}{\p^2}(a+ \p a)$.\\
 In this situation, $\theta(P(z)) = (x'(P(z), y'(P(z))$ satisfies  $x'(P(z))=0$ and therefore we have $\delta_{\psi_{a,b},K_\p} \circ \theta (P(z)) = ( \pm 2 \p^3d\sqrt{a}, \  \pm \frac{1}{2 \p^3d\sqrt{a}})$.
 Now note that
 \begin{align*}
     \p^{-3} (\pm 2 \p^3 d\sqrt{a}) & = \pm 2  d\sqrt{a} \equiv \pm 8 a\sqrt{a} \pmod{\p^3} = (\pm 2 \sqrt{a})^3 \pmod{\p^3}.
 \end{align*}
Therefore, in this case, we obtain $\overline{\p^{-\upsilon_\p(C)}(\pm 2\sqrt{a})^{-3}C}=\overline{1} \in \big({\OO_{K_\p}^*}/{U_{K_\p}^3}\big)^{3}$. \qed
\end{enumerate}
\end{enumerate}

 \subsubsection{\underline{{\bf Case-2: $3 \mid b$} }}\label{case2type2}
In this case, we assume that $3 \mid b$ (this implies $3 \nmid a$) and study the relationship between the groups $\frac{\widehat{E}_{a,b}(K_\p)}{\psi_{a,b}(E_{a,b}(K_\p))}$, 
$\frac{V_3}{A_\p^{*3}}$ and $\Big(\frac{A_\p^*}{A_\p^{*3}}\Big)_{N=1}$.
 
 By Tate's algorithm, $E_{a,b}:y^2=x^3+a(x-b)^2$ is $\p$-minimal and has multiplicative reduction of type $I_n$ at $\p$, where $n= 3 \upsilon_\p(b)$. Thus, the Tamagawa number 
 $c_\p(E_{a,b}) = \begin{cases}
  2,  & \text{ if } a \notin K_\p^{*2}, \\
  3 \ \upsilon_\p(b), & \text{ if } a \in K_\p^{*2}.
\end{cases}$

Observe that the equation of  $\widehat{E}_{a,b}:y^2=x^3-27a(x-d)^2$ is not $\p$-minimal and by Tate's algorithm, the $\p$-minimal model for $\widehat{E}_{a,b}$ is given by $E' : Y^2= X^3+ a_2 X^2 + b_2 X + c_2$, where $a_2$, $b_2$, $c_2$ are as defined in \eqref{eq:defnofE2}. As in the previous case,  $\theta: E' \to \widehat{E}_{a,b}$ as in \eqref{eq:defoftheta} defines a $K$-isomorphism and the isogeny $\Psi_{a,b} : E_{a,b} \to E'$ is again given by \eqref{eq:newpsi}.

Moreover, $E'$ has multiplicative reduction of type $I_n$ at $\p$, where $n=  \upsilon_\p(b)$. Thus,
 $c_\p(E') = 2$, if $a \notin K_\p^{*2}$ and $c_\p(E')=\upsilon_\p(b)$, if $a \in K_\p^{*2}$.
In particular, note that $c_\p(E_{a,b})=3c_\p(\widehat{E}_{a,b})$, when $a \in K_\p^{*2}$. Using formal groups, we have that $\Psi'_{a,b}(0)=-{\zeta_3}$ and hence, $||\Psi'_{a,b}(0)||_\p^{-1}=1$.

\begin{lemma}\label{3midb}
For $3 \mid b$, we have
(whether or not $a \in K_\p^{*2}$) that  $\Big|\frac{\widehat{E}_{a,b}(K_\p)}{\psi_{a,b}(E_{a,b}(K_\p))}\Big| =1$ and thus trivially, $\frac{\widehat{E}_{a,b}(K_\p)}{\psi_{a,b}(E_{a,b}(K_\p))} \subsetneq \Big(\frac{A_\p^*}{A_\p^{*3}}\Big)_{N=1}$. 
\end{lemma}
\begin{proof}
  This follows easily from the above discussion and using \eqref{eq:schaeferfor}.
\end{proof}

\subsubsection{\underline{{\bf Case-3: $3 \mid a$}}}\label{case3type2}
We finally consider the case $3 \mid a$;  this implies $3 \nmid b$. We study the bounds for $\frac{\widehat{E}_{a,b}(K_\p)}{\psi_{a,b}(E_{a,b}(K_\p))}$ using $\frac{V_3}{A_\p^{*3}}$ and $\Big(\frac{A_\p^*}{A_\p^{*3}}\Big)_{N=1}$.

\noindent If $a \notin K_\p^{*2}$, then we know that $\{1\}=\frac{V_3}{A_\p^{*3}} \subset \frac{\widehat{E}_{a,b}(K_\p)}{{\psi_{a,b}}(E_{a,b}(K_\p))} \subset \Big(\frac{A_\p^*}{A_\p^{*3}}\Big)_{N=1}$ from Props. \ref{ais-1modulo3} and  \ref{sizeofV3}.

From now on, we take $a \in K_\p^{*2}$. First, we consider the special case when $3^6 \mid d = 4a+27b$;  then $27 \mid \mid a$. Setting $a_1 = - \frac{a}{27}$ and $b_1 = \frac{d}{3^6}$, we see that $\widehat{E}_{a,b} \cong {E}_{a_1,b_1}$ and ${E}_{a,b} = \widehat{E}_ {a_1,b_1}$. Note here that $3 \nmid a_1$ and $a_1 \in K_\p^{*2}$. Then following a similar proof as in Case 1 \S\ref{case1type2} and Case 2 \S\ref{case2type2}, it follows that:

\noindent (1) If $3^6 \mid \mid d$ or equivalently, $3 \nmid b_1$, we see that 
$\frac{V_3}{A_\p^{*3}} \cong \frac{\widehat{E}_{a,b}(K_\p)}{{\psi}_{a,b}(E_{a,b}(K_\p))} \subset \Big(\frac{A_\p^*}{A_\p^{*3}}\Big)_{N=1}$.

\noindent (2) If $3^7 \mid d$ or equivalently $3 \mid b_1$, then $\{1\}=  \frac{\widehat{E}_{a,b}(K_\p)}{{\psi_{a,b}}(E_{a,b}(K_\p))} \subsetneq \Big(\frac{A_\p^*}{A_\p^{*3}}\Big)_{N=1}$.

\vspace{2mm}
We summarize this in the following remark:
\begin{rem}\label{remark6.8}
For $a \notin K_\p^{*2}$, we have $|V_3/A_\p^{*3}|=1$ by Prop. \ref{sizeofV3}.  Then
\begin{enumerate}
    \item If $3 \nmid a$, then $\frac{V_3}{A_\p^{*3}} = \delta_{{\psi}_{a,b},K_\p}(\widehat{E}_{a,b}(K_\p)) \subsetneq \Big(\frac{A_\p^*}{A_\p^{*3}}\Big)_{N=1}$ by Props. \ref{3nmidab1}(1) and \ref{3midb}(1).
    \item If $3 \mid a$, then $\frac{V_3}{A_\p^{*3}} \subset \delta_{\psi_{a,b},K_\p}(\widehat{E}_{a,b}(K_\p)) \subset \Big(\frac{A_\p^*}{A_\p^{*3}}\Big)_{N=1}$ by Prop. \ref{ais-1modulo3}.
\end{enumerate}
For $a \in K_\p^{*2}$, we have $|V_3/A_\p^{*3}|=3$ by Prop. \ref{sizeofV3}. Then
\begin{enumerate}
    \item If $3 \nmid ab$, then $\frac{V_3}{A_\p^{*3}}=\delta_{{\psi}_{a,b},K_\p}(\widehat{E}_{a,b}(K_\p)) \subsetneq \Big(\frac{A_\p^*}{A_\p^{*3}}\Big)_{N=1}$ by Prop. \ref{3nmidab2}.
    \item If $3 \mid b$, then $\{1\}=\delta_{{\psi}_{a,b},K_\p}(\widehat{E}_{a,b}(K_\p))  \subsetneq \Big(\frac{A_\p^*}{A_\p^{*3}}\Big)_{N=1}$ by Prop. \ref{3midb}(2).
    \item If $27 \mid\mid a$ and $3^6 \mid\mid d$, then $\frac{V_3}{A_\p^{*3}} \cong \delta_{{\psi}_{a,b},K_\p}(\widehat{E}_{a,b}(K_\p)) \subsetneq \Big(\frac{A_\p^*}{A_\p^{*3}}\Big)_{N=1}$ by Props. \ref{3nmidab1}(2) and \ref{3nmidab2}.
    \item If $27 \mid\mid a$ and $3^7 \mid d$, then $\{1\}= \delta_{{\psi}_{a,b},K_\p}(\widehat{E}_{a,b}(K_\p))  \subsetneq \Big(\frac{A_\p^*}{A_\p^{*3}}\Big)_{N=1}$ by Lemma \ref{3midb}(2).
\end{enumerate} 
\end{rem}

\section{Bounds for the $\psi_{a,b}$- and $3$-Selmer Groups of $E_{a,b}/K$}\label{type2global}
In this subsection, we give bounds on $\dim_{\F_3}{\rm Sel}^{\psi_{a,b}}(E_{a,b}/K)$ and $\dim_{\F_3}{\rm Sel}^3(E_{a,b}/K)$ in terms of $\dim_{\F_3}Cl_S(L)[3]$ for certain sets $S$ that we define shortly.

\begin{rem}\label{rmkforboundsonpsihat}
As explained in \S\ref{iso}, we have ${\rm Sel}^{\widehat{\psi}_{a,b}}(\widehat{E}_{a,b}/K) \cong {\rm Sel}^{{\psi}_{-27a,d}}(E_{-27a,d}/K)$ using $\widehat{\psi}_{a,b}=\lambda_{27} \circ \psi_{-27a,d}$. Further, writing $gcd(-27a,d)=t^2g$, where $g$ is square-free in $K$, $\widehat{E}_{a,b}$ can be identified with $E_{{-27a}/{t^2}, {d}/{t^2}}$ upto an isomorphism defined over $K$.
Thus to compute bounds of $\dim_{\F_3} {\rm Sel}^{\widehat{\psi}_{a,b}}(\widehat{E}_{a,b}/K)$, it suffices to compute the same for ${\rm Sel}^{\psi_{a,b}}(E_{a,b}/K)$ for all $a,b \in \Z$ with $gcd(a,b)$ square-free in $K$.
\end{rem}

 Recall that $a,b \in \Z$ with $gcd(a,b)$ square-free in $K$;  In particular, $3 \nmid gcd(a,b)$. Set 
$$T_{2,1}:= \begin{cases} \{2\OO_K\}, & \text{ if }  a \in K_2^{*2}, 2 \nmid a \text{ and }  4 \nmid b, \\ \emptyset,& \text{ otherwise. }  \end{cases} \quad  
T_{3,1}:=\begin{cases} \{\p\}, & \text{ if } a \in K_\p^{*2},  \p \mid a \text{ and } \upsilon_\p(a) \neq 6,\\
\{\p\}, & \text{ if } a \in K_\p^{*2},  \upsilon_\p(a) = 6 \text{ and } \upsilon_\p(d) < 12,\\
\emptyset, & \text{ otherwise.} \end{cases}$$
$$T_{2,2}:= \begin{cases} \{2\OO_K\}, & \text{ if }  a \in K_2^{*2}, 2 \nmid a \text{ and } 8 \mid b, \\ \emptyset,& \text{ otherwise. }  \end{cases} \quad \enspace 
T_{3,2}:= \begin{cases} \{\p\}, & \text{ if }  a \in K_\p^{*2}, \p \mid b \text{ and } \p \nmid a, \\ \emptyset,& \text{ otherwise. }  \end{cases}$$
$$T_{2,3}:= \begin{cases} \{2\OO_K\}, & \text{ if }  a \in K_2^{*2}, 2 \nmid a \text{ and } 4 \mid\mid b, \\ \emptyset,& \text{ otherwise. }  \end{cases} \quad  \enspace 
T_{3,3}:= \begin{cases} \{\p\}, & \text{ if }  a \in K_\p^{*2} \text{ and } \upsilon_\p (d) > 12, \\ \emptyset,& \text{ otherwise. }  \end{cases}$$

\begin{defn}\label{defofS1S2S3}
We define the following subsets of $\Sigma_K$, the set of all finite primes of $K$:
$$S_1= \Big\{ \q \in \Sigma_K \setminus \{ \p \} \mid  a \in K_\q^{*2}, \enspace \q \mid a \ \text{ and } \ \upsilon_\q(4ab^2) \not\equiv 0 \pmod{6} \Big\} \ \cup T_{2,1} \ \cup T_{3,1},$$
$$S_2= \Big\{ \q \in \Sigma_K \setminus \{\p\} \mid  a \in K_\q^{*2}, \enspace \q \mid b \ \text{ and } \ \q \nmid 2a  \Big\} \ \cup T_{2,2} \ \cup T_{3,2} \quad \text{ and }$$
$$S_3= \Big\{ \q \in \Sigma_K \setminus \{ \p \} \mid  a \in K_\q^{*2}, \enspace \q \mid d \ \text{ and } \ \q \nmid 2a \Big\} \ \cup T_{2,3} \ \cup T_{3,3}.$$
\end{defn}
Observe that $S_1, \ S_2$ and $S_3$ are pairwise disjoint finite sets of primes of $K$.
\begin{defn}\label{defofSLtype2}
Let $L$ be a field. Define the following subsets of $\Sigma_L$:
$$S_{1,2}(L):=\{\mathfrak{Q} \in \Sigma_L \mid \mathfrak{Q} \cap \OO_K \in S_1 \cup S_2 \} \ \text{ and } \ 
S_{1,3}(L):=\{\mathfrak{Q} \in \Sigma_L \mid \mathfrak{Q} \cap \OO_K \in S_1 \cup S_3 \}.$$
\end{defn}

\begin{rem}
Using Tate's algorithm, we observe that if $\q \in S_1 \cup S_2$ (resp. $\q \in S_1 \cup S_3$) including $\q=\p$, then $E_{a,b}$ and $\widehat{E}_{a,b}$ have either split multiplicative reduction or additive reduction of Kodaira type \rom{4}  or \rom{4}$^*$ at $\q$. 
\end{rem}
\noindent The $\F_3$-subspaces $M(S_1 \cup S_2, a)$,  $N(S_1 \cup S_3, a)$ and $N'(S_1 \cup S_3, a)$ of $\big(L^*/L^{*3}\big)_{N=1}$ are defined in Prop. \ref{newM1M2sel}. 
\begin{theorem}\label{containment2}
We have the inclusion of $\F_3$-modules and the inequalities of their $3$-ranks:
\begin{enumerate}
    \item If $a \notin K^{*2}$, then $M(S_1 \cup S_2, a) \subset {\rm Sel}^{\psi_{a,b}}(E_{a,b}/K) \subset N(S_1 \cup S_3, a)$.  Hence, we get that $$h^3_{S_{1,2}(L)} \le \dim_{\F_3} {\rm Sel}^{\psi_{a,b}}(E_{a,b}/K) \le h^3_{S_{1,3}(L)} + |S_{1,3}(L)| +2.$$
    \item If $a \in K^{*2}$, then ${\rm Sel}^{\psi_{a,b}}(E_{a,b}/K) \subset N'(S_1 \cup S_3, a)$. Hence, $\dim_{\F_3} {\rm Sel}^{\psi_{a,b}}(E_{a,b}/K)  \le |S_1 \cup S_3| +1.$ 
\end{enumerate}
\end{theorem}

\begin{proof}
We start with the case $a \notin K^{*2}$. From Props. \ref{type2nsq}, \ref{type2redn} and Table \ref{tab:type2atq} (for $\q \nmid 3$) and Props. \ref{ais-1modulo3}, \ref{3nmidab1}, \ref{3midb} and Remark \ref{remark6.8} 
(for $\p \mid 3$), it follows that $\text{Im } \delta_{{\psi}_{a,b},K_\q} \subset (A_\q^*/A_\q^{*3})_{N=1}$ for all $\q \notin S_1 \cup S_3$. Thus, by Prop. \ref{newM1M2sel}, ${\rm Sel}^{\psi_{a,b}}(E_{a,b}/K) \subset N(S_1 \cup S_3, a)$. 

On the other hand, it follows from Props. \ref{type2nsq}, \ref{type2redn} and Table \ref{tab:type2atq}, that $(A_\q^*/A_\q^{*3})_{N=1} \subset \text{Im } \delta_{{\psi}_{a,b}, K_\q}$ for $\q \notin S_1 \cup S_2 \cup \{\p\}$. Further, from Props. \ref{3nmidab1}, \ref{3midb}, Remark \ref{remark6.8}, 
it follows that $V_3/A_\p^{*3} \subset \text{Im } \delta_{{\psi}_{a,b}, K_\p}$, if $\p \notin S_1 \cup S_2$. Thus, from Prop. \ref{newM1M2sel}, we get $M(S_1 \cup S_2, a) \subset {\rm Sel}^{\psi_{a,b}}(E_{a,b}/K)$.

Now if $a \in K^{*2}$, then from Prop. \ref{type2redn} and Table \ref{tab:type2atq} (for $\q \nmid 3$) and Props. \ref{3nmidab1}, \ref{3midb} and Remark \ref{remark6.8} 
(for $\p \mid 3$), it follows that $\text{Im }(\delta_{{\psi}_{a,b},K_\q}) \subset (A_\q^*/A_\q^{*3})_{N=1}$ for all $\q \notin S_1 \cup S_3$. Hence, by Prop. \ref{newM1M2sel}, ${\rm Sel}^{\psi_{a,b}}(E_{a,b}/K) \subset N'(S_1 \cup S_3, a)$. 

\noindent In both cases (1) and (2), the statements about ranks follow directly from Prop. \ref{boundsgen}.
\end{proof}

\begin{rem}\label{rmktopsihat}
Following Remark \ref{rmkforboundsonpsihat}, the bounds for ${\rm Sel}^{\widehat{\psi}_{a,b}}(\widehat{E}_{a,b}/K)$ are computed similarly: 
\begin{enumerate}
	\item If $a \notin K^{*2}$, then $h^3_{S_{1,3}(L)} \le \dim_{\F_3} {\rm Sel}^{\widehat{\psi}_{a,b}}(\widehat{E}_{a,b}/K)  \le h^3_{S_{1,2}(L)}+ |S_{1,2}(L)|+2.$
	\item If $a \in K^{*2}$, then  $\dim_{\F_3} {\rm Sel}^{\widehat{\psi}_{a,b}}(\widehat{E}_{a,b}/K)  \le |S_1 \cup S_2| +1.$ 
\end{enumerate}
\end{rem}
 \noindent Theorem \ref{containment2} and Remark \ref{rmktopsihat} have the following consequence on ${\rm Sel}^3(E_{a,b}/K)$:
\begin{corollary}\label{type2sel3}
If $a \notin K^{*2}$, then from \cite[Lemma~6.1]{ss} we have the exact sequence 
\begin{equation}\label{eq:exactseqfortype2}
    0 \to {\rm Sel}^{{\psi}_{a,b}}({E}_{a,b}/K) \to {\rm Sel}^{3}({E}_{a,b}/K) \to {\rm Sel}^{{\widehat{\psi}}_{a,b}}(\widehat{E}_{a,b}/K) \to \frac{\Sh(\widehat{E}_{a,b}/K)[\widehat{\psi}_{a,b}]}{\psi_{a,b}(\Sh(E_{a,b}/K)[3]} \to 0
\end{equation} 
and it follows that 
$$\text{max } \{ h^3_{S_{1,2}(L)}, \text{rk } E_{a,b}(K) \} \leq \dim_{\F_3} {\rm Sel}^3(E_{a,b}/K) \leq h^3_{S_{1,2}(L)}+h^3_{S_{1,3}(L)}+|S_{1,2}(L)|+|S_{1,3}(L)|+4.$$
If $a \in K^{*2}$, then $\text{rk } E_{a,b}(K) \le \dim_{\F_3} {\rm Sel}^{3}(E_{a,b}/K) \le |S_1 \cup S_2|+|S_1 \cup S_3|+2$.
\qed
\end{corollary}

\begin{theorem}[Refined Bounds]\label{cortotype2}
If $a \notin K^{*2}$ and $3 \nmid a$, then we have that 
\begin{multline*}
    \text{max } \{h^3_{S_{1,2}(L)}, \enspace h^3_{S_{1,3}(L)} +|S_3|-|S_2|-1\} \leq \dim_{\F_3} {\rm Sel}^{\psi_{a,b}}(E_{a,b}/K) \\
    \leq \text{min } \{h^3_{S_{1,2}(L)} +|S_{1,2}(L)|+|S_3|-|S_2|+1, \enspace h^3_{S_{1,3}(L)}+|S_{1,3}(L)|+2 \}.
\end{multline*}
\noindent In particular, if $S_1 \cup S_2 \cup S_3$ is empty, then $\dim_{\F_3} {\rm Sel}^{\psi_{a,b}}(E_{a,b}/K) \in \{h^3_L, \enspace h^3_L+1\}$.
\end{theorem}

\proof
Using \cite[Eq.~(1.22), (1.26) and (3.4)]{cass2}, we get that
\begin{equation}\label{eq:casselsfor}
\frac{|{\rm Sel}^{\widehat{\psi}_{a,b}}(\widehat{E}_{a,b}/K)|}{|{\rm Sel}^{\psi_{a,b}}(E_{a,b}/K)|}=\frac{|\widehat{E}_{a,b}(K)[\widehat{\psi}_{a,b}]|}{|E_{a,b}(K)[\psi_{a,b}]|} \cdot \frac{|E_{a,b}(\mathbb{C})[\psi_{a,b}]|}{|\widehat{E}_{a,b}(\mathbb{C})/\psi_{a,b}(E_{a,b}(\mathbb{C}))|} \cdot \prod_{\q \in \Sigma_K}{\frac{|E_{a,b}(K_\q)[\psi_{a,b}]|}{|\widehat{E}_{a,b}(K_\q)/\psi_{a,b}(E_{a,b}(K_\q))|}}.
\end{equation}
\noindent Also $\big|\frac{\widehat{E}_{a,b}(K_\q)}{\psi_{a,b}(E_{a,b}(K_\q))}\big|=\frac{||\psi_{a,b}'(0)||_\q^{-1} \cdot |E_{a,b}(K_\q)[\psi_{a,b}]|\cdot c_\q(\widehat{E}_{a,b})}{c_\q(E_{a,b})}$, where 
$||\psi_{a,b}'(0)||_\q^{-1}=1$ for any $\q \in \Sigma_K$ by \eqref{eq:schaeferfor}.
Now, note that $|{E}_{a,b}(\mathbb{C})[{\psi}_{a,b}]|=3$ and $\psi_{a,b}:E_{a,b}(\mathbb{C}) \to \widehat{E}_{a,b}(\mathbb{C})$ is surjective.
Thus, \eqref{eq:casselsfor} becomes 
$\frac{|{\rm Sel}^{\widehat{\psi}_{a,b}}(\widehat{E}_{a,b}/K)|}{|{\rm Sel}^{\psi_{a,b}}(E_{a,b}/K)|}=3 \cdot \underset{\q \in \Sigma_K}{\prod} \frac{c_\q(E_{a,b})}{c_\q(\widehat{E}_{a,b})}$.

Define $S':=\{ \q \in \Sigma_K \mid c_\q(E_{a,b})=3 \ c_\q(\widehat{E}_{a,b}) \}$ and $S'':=\{ \q \in \Sigma_K \mid c_\q(\widehat{E}_{a,b}) = 3 \ c_\q(E_{a,b}) \}$. 
Then we have $\frac{|{\rm Sel}^{\widehat{\psi}_{a,b}}(\widehat{E}_{a,b}/K)|}{|{\rm Sel}^{\psi_{a,b}}(E_{a,b}/K)|}=3^{|S'|-|S''|+1}$.
Under the assumption $3 \nmid a$, it is easy to see from the local theory of the curves $E_{a,b}$ (Table \ref{tab:type2atq}, \S\ref{case2type2}) that $S'=S_2$ and $S''=S_3$. 
This implies $\frac{|{\rm Sel}^{\widehat{\psi}_{a,b}}(\widehat{E}_{a,b}/K)|}{|{\rm Sel}^{\psi_{a,b}}(E_{a,b}/K)|}=3^{|S_2|-|S_3|+1}$. Rearranging the terms and using Remark \ref{rmktopsihat} we get that 
$$\dim_{\F_3} {\rm Sel}^{\psi_{a,b}}(E_{a,b}/K)= \dim_{\F_3} {\rm Sel}^{\widehat{\psi}_{a,b}}(\widehat{E}_{a,b}/K) + |S_3| - |S_2| -1 \geq h^3_{S_{1,3}(L)} +|S_3|-|S_2|-1.$$ 
This along with Theorem \ref{containment2} gives that $\dim_{\F_3} {\rm Sel}^{\psi_{a,b}}(E_{a,b}/K) \geq \text{max}\{h^3_{S_{1,2}(L)}, \enspace h^3_{S_{1,3}(L)} +|S_3|-|S_2|-1\}$.

Similarly by Remark \ref{rmktopsihat}, we have that $$\dim_{\F_3} {\rm Sel}^{\widehat{\psi}_{a,b}}(\widehat{E}_{a,b}/K) + |S_3| - |S_2| -1 = \dim_{\F_3} {\rm Sel}^{{\psi}_{a,b}}({E}_{a,b}/K)  \le h^3_{S_{1,2}(L)}+ |S_{1,2}(L)|+ |S_3| - |S_2| +1.$$ 
Together with Theorem \ref{containment2}, this gives that $$\dim_{\F_3} {\rm Sel}^{\psi_{a,b}}(E_{a,b}/K) \leq \text{min} \{h^3_{S_{1,2}(L)} +|S_{1,2}(L)|+|S_3|-|S_2|+1, \enspace h^3_{S_{1,3}(L)}+|S_{1,3}(L)|+2 \}.$$
This completes the proof. \qed

\vspace{3mm}
\noindent {\bf Example (1):} Take $a=29$ and $b=76$. Then $S_1=S_2=\emptyset$ and $S_3=\{2\OO_K\}$;  also, $h^3_{S_{1,2}(L)}=h^3_L=1$ and $h^3_{S_{1,3}(L)}=0$. So, by Theorem \ref{cortotype2}, $1 \le \dim_{\F_3}{\rm Sel}^{\psi_{a,b}}(E_{a,b}/K) \le 3$.

\vspace{1mm}
\noindent {\bf Example (2):} For $a=79$ and $b=79$, observe that $S_1=S_2=S_3=\emptyset$ and $h^3_{S_{1,2}(L)}=h^3_{S_{1,3}(L)}=h^3_L=2$. So, by Theorem \ref{cortotype2}, $2 \le \dim_{\F_3}{\rm Sel}^{\psi_{a,b}}(E_{a,b}/K) \le 3$.

\vspace{1mm}
\noindent {\bf Example (3):} Taking $a=2263$ and $b=72$, we have that $S_1=S_2=S_3=\emptyset$ and $h^3_L=3$. Thus, we get that $3 \le \dim_{\F_3}{\rm Sel}^{\psi_{a,b}}(E_{a,b}/K) \le 4$ using Theorem \ref{cortotype2}.


\vspace{1mm}
\noindent {\bf Example (4):} If $a=113$ and $b=22$, then $S_1=\{2\OO_K\}$, $S_2=\{11\OO_K\}$ and $S_3=\emptyset$. We get that $Cl_{S_{1,2}(L)(L)} =\Z/6$ and $Cl_{S_{1,3}(L)(L)} =\Z/3$. Thus, Theorem \ref{cortotype2} gives that $1 \le \dim_{\F_3}{\rm Sel}^{\psi_{a,b}}(E_{a,b}/K) \le 5$.

\begin{corollary}\label{3.9}
    Assume that $a \notin K^{*2}$ and $3 \nmid a$. Then the global root number $\omega(E_{a,b}/K)$ is determined by the parity of $|S_2| + |S_3|$. More precisely,
    \begin{equation}
        \omega(E_{a,b}/K) = (-1)^{|S_2|+|S_3|+1} = \begin{cases}
            -1 , & \text{ if } |S_2| \equiv |S_3|  \pmod{2}, \\
            1 , & \text{ if } |S_2| \not\equiv |S_3|  \pmod{2}. \\ 
        \end{cases}
    \end{equation}
\end{corollary}

\begin{proof}
From the exact sequence \eqref{eq:exactseqfortype2}, we obtain 
\begin{equation}\label{eq:simplifiedexactseq}
\dim_{\F_3} {\rm Sel}^{3}({E}_{a,b}/K) \equiv \dim_{\F_3} {\rm Sel}^{{\psi}_{a,b}}({E}_{a,b}/K) + \dim_{\F_3} {\rm Sel}^{{\widehat{\psi}}_{a,b}}(\widehat{E}_{a,b}/K) \pmod{2}
\end{equation}
since $\dim_{\F_3} \frac{\Sh(\widehat{E}_{a,b}/K)[\widehat{\psi}_{a,b}]}{\psi_{a,b}(\Sh(E_{a,b}/K)[3]}$ is even \cite[Prop.~49]{bes}. Also, from the proof of Theorem \ref{cortotype2}, we get 
$$\dim_{\F_3} {\rm Sel}^{\psi_{a,b}}(E_{a,b}/K)= \dim_{\F_3} {\rm Sel}^{\widehat{\psi}_{a,b}}(\widehat{E}_{a,b}/K) + |S_3| - |S_2| -1$$
and hence $\dim_{\F_3} {\rm Sel}^{3}({E}_{a,b}/K) \equiv |S_3| - |S_2| -1 \pmod{2}$. Now $E_{a,b}(K)[3]=0$ since $a \notin K^{*2}$ by assumption. It follows that $r= \text{corank } {\rm Sel}^{3^\infty}({E}_{a,b}/K) \equiv \dim_{\F_3} {\rm Sel}^{3}({E}_{a,b}/K) \pmod{2}$. The result follows from $3$-parity conjecture of $E_{a,b}$ over $K$ which is known due to Nekov{\'a}{\v r} and Dokchitser-Dokchitser.
\end{proof}

By studying the image of the Kummer maps over $\Q_\ell$, following a method similar to that in Section \ref{type2curves}, one can obtain an upper bound of $\dim_{\F_3} {\rm Sel}^{\psi_{a,b}}(E_{a,b}/\Q)$. In particular, if this upper bound is small, it is possible to determine the dimension of the $3$-Selmer group of $E_{a,b}$ over $\Q$. Using this approach, we provide an infinite family of elliptic curves $E_{a,b}$ with $3$-Selmer rank equal to $1$.

\begin{proposition}\label{3.10}
   For a positive proportion of natural numbers $n$, one has that
   $$\omega(E_{n,n}/\Q)=-1, \ \omega(E_{-3n,-3n}/\Q)=+1 \text{ and } \dim_{\F_3} {\rm Sel}^3(E_{n,n}/\Q) =1.$$
\end{proposition}

\begin{proof}
    Let $n$ be a square-free natural number satisfying $n \equiv 2 \text{ or } 3 \pmod{4}$, $n \equiv 2 \pmod 3$  and $n$ is a quadratic non-residue modulo $31$. Let $K_n:=\Q(\sqrt{-3n})$ and $\OO_{K_n}$ be its ring of integers. Consider the elliptic curve $E_{n,n}$. Then for the above choice of $n$, for every prime $\ell$, using \eqref{eq:schaeferfor} one  can show that
    $\text{Im } \delta_{{\psi}_{n,n},\Q_\ell} \subset (A_\ell^*/A_\ell^{*3})_{N=1}$, where $A_{\ell}:= \frac{\Z_\ell[X]}{(X^2+3n)}$. Consequently, we have that $\dim_{\F_3} {\rm Sel}^{\psi_{n,n}}(E_{n,n}/\Q)$ is bounded by $h^3_{K_n}+ \dim_{\F_3} \frac{\mathcal{O}_{K_n}^*}{\mathcal{O}_{K_n}^{*3}} = h^3_{K_n}$.
    As $|E_{n,n}(\R)[\psi_{n,n}]|=3$, from Cassels' formula \eqref{eq:casselsfor} and Schaefer's formula \eqref{eq:schaeferfor}, we obtain that $\frac{|{\rm Sel}^{\widehat{\psi}_{n,n}}(\widehat{E}_{n,n}/\Q)|}{|{\rm Sel}^{\psi_{n,n}}(E_{n,n}/\Q)|}=3$. Using this, we deduce from  \cite[Lemma 6.1]{ss} (see also \eqref{eq:simplifiedexactseq})  that $\dim_{\F_3} {\rm Sel}^3(E_{n,n}/\Q) \equiv 1 \pmod 2$ and as $E_{n,n}(\Q)[3]=0$, we have that $\text{corank } {\rm Sel}^{3^\infty}({E}_{n,n}/\Q) \equiv \dim_{\F_3} {\rm Sel}^3(E_{n,n}/\Q) \equiv 1 \pmod 2$. It follows that $\omega(E_{n,n}/\Q) = -1$, by the $3$-parity conjecture. Further, observe that $E_{n,n}$ is isomorphic to $E_{-3n,-3n}$ over $K$ and $\omega(E_{n,n}/K) =-1$ by Corollary \ref{3.9}. Hence, we get that $\omega(E_{-3n,-3n}/\Q)= +1$.
    
    By a result of Nakagawa-Horie (which can be conveniently found in \cite[Lemma 2.2]{by}), we see that in each of the $30$ equivalent classes, $n \equiv a \pmod{372 }$ determined by the conditions  $n \equiv 2 \text{ or } 3 \pmod{4}$, $n \equiv 2 \pmod 3$ and $n$ is a quadratic non-residue modulo $31$, $h^3_{K_n}$ vanishes for at least $50\%$ of $n$. Further, if $h^3_{K_n}=0$, then by the discussion above, we get that ${\rm Sel}^{\psi_{n,n}}(E_{n,n}/\Q)$ is trivial and $\dim_{\F_3} {\rm Sel}^3(E_{n,n}/\Q) =1$. 

    Let $m, M \in \N$ be such that $d:=\text{gcd}(m,M)$ is square free. Then we have \begin{small}$$\sum_{\substack{0< n \le X \\ n \equiv m \pmod{M}}} \mu^2(n) = \frac{X}{M \zeta(2)} \Big(\prod_{p \mid \frac{M}{d}} \frac{1}{1- p^{-2} } \Big) + o(X).  $$\end{small} It follows that the number of  square-free integers \begin{small}$0 < n  \le X$,\end{small} with \begin{small}$n \equiv 2, 11 \pmod{12}$\end{small} and $n$ is not a square modulo $31$ is of the order of $\frac{30X}{372 \zeta(2)} \cdot \frac{2^2}{2^2-1} \cdot \frac{3^2}{3^2-1}\cdot \frac{31^2}{31^2-1} = \frac{31X}{2^8 \zeta(2)}$.
    
    Putting all of it together, we see that for $X \gg 0$, the number of square-free natural numbers $n \le X$, such that $\omega(E_{n,n}/\Q)=-1$, $\omega(E_{-3n,-3n}/\Q)=+1$ and $\dim_{\F_3} {\rm Sel}^3(E_{n,n}/\Q) =1$ is at least $\frac{31X}{2^9 \zeta(2)}$. This completes the proof.
\end{proof}

\section{Arbitrary large rank of $3$-Selmer groups of elliptic curves over $K$}\label{largeselmerrank}
We construct  infinite families of elliptic curves with  arbitrary large $3$-Selmer group over $K$. 
\begin{theorem}\label{large3selmer}
Given any $n \geq 0$, there exist infinitely many elliptic curves $E_{a,1}:y^2=x^3+a(x-1)^2$ with $a \in \Z$ such that $\dim_{\F_3} {\rm Sel}^3(E_{a,1}/K) \geq \dim_{\F_3} {\rm Sel}^{\psi_{a,1}}(E_{a,1}/K) \geq 2n$. These curves have a rational $3$-isogeny and no non-trivial $3$-torsion points over $K$.  
\end{theorem}

\begin{proof}
Take $a=\frac{p_1 p_2 \cdots p_{2n+1}-27}{4}$, where $p_i$'s are all distinct rational primes and each $p_i \equiv -1 \pmod{12}$. Then $p_1 p_2 \cdots p_{2n+1} = 12s -1$ for some $s$ and  $a = \frac{12s-28}{4} = 3s -7 \equiv -1 \pmod 3$. Hence, $a \notin K^{*2}$;  moreover,  $3 \nmid a$ and $a \notin K_\p^{*2}$. Consider the sets $S_2$, $S_3$ from the Definition \ref{defofS1S2S3} and observe that $S_2 = \emptyset$. Further, $d:=4a+27=p_1 p_2 \cdots p_{2n+1}$, so  each $p_i \mid d$ and $p_i \nmid 2a$. Also, since each $p_i \equiv -1 \pmod{12}$, we see that each $p_i$ is inert in $K$ and $\big(\frac{-3}{p_i}\big)=-1$ for all $i$. Hence, either $\big(\frac{a}{p_i}\big)=1$ or $\big(\frac{-3a}{p_i}\big)=1$, which implies  $a \in K_{p_i}^{*2}$ for all $p_i$. This in turn implies that each $p_i \in S_3$;  in fact, $S_3=\{p_1, \cdots, p_{2n+1}\}$. Indeed, if $\q \neq p_i$ for any $i$, then $\q \nmid d$, so $\q \notin S_3$. As $a \notin K_\p^{*2}$, \ $\p \notin S_3$ and $2\OO_K \notin S_3$. Thus, we have $|S_3|=2n+1$. By Theorem \ref{cortotype2}, we get that $\dim_{\F_3} {\rm Sel}^{\psi_{a,1}}(E_{a,1}/K) \geq 2n$.

Clearly, there are infinitely many choices of $a$ for the given $n$, hence the statement holds.
\end{proof}

\begin{rem}[Related works]\label{largecomp}
In \cite{cass3}, Cassels showed that $\dim_{\F_3}{\rm Sel}^3(E/\Q)$ for a class of elliptic curves $E$, each having a rational $3$-torsion point and $j$-invariant $0$, can be arbitrarily large.
For $p=2$, \cite{bo, kr} prove that for any given $n>0$, there exists an elliptic curve $E$ such that $\dim_{\F_2}{\rm Sel}^2(E/\Q) \ge n$, assuming that $E[2] \subset E(\Q)$.
On the other hand, for certain $n \geq 5$, there are several results (see \cite{ks, fish}) proving that the $n$-Selmer groups of elliptic curves over $\Q$ can be arbitrarily large. 

In contrast, in our result on $\dim_{\F_3}{\rm Sel}^3(E_{a,b}/K)$ in Theorem \ref{large3selmer}, the elliptic curves $E_{a,b}/\Q$ have a $\Q$-rational $3$-isogeny, do not have a non-trivial $3$-torsion point over $K$ and the $j$-invariant $j(E_{a,b}) \neq 0$.  

In their recent article \cite[Theorem~1.3]{bkos}, the authors show that for a fixed number field $F$ and $r \geq 1$, almost all elliptic curves $E/F$ with a cyclic $F$-rational $9$-isogeny have a positive proportion of quadratic twists $E_s$ with $\Sh(E_s/F)[3] \geq 9^r$ (see \cite[Theorem~1.3]{bkos} for a  precise statement). 
In a different direction, given a prime $p$ and an elliptic curve $E$ defined over a number field $F$, it was shown by \cite{ces} that $\dim_{\F_p}{\rm Sel}^p(E/F')$ is unbounded where $F'$ varies over $\Z/p\Z$-extensions of $F$.
There are important works of Mazur-Rubin on large corank of  $p^\infty$-Selmer groups.
\end{rem}

We have the following interesting application of Theorem \ref{containment2}.
\begin{corollary}\label{biquadfld}
Given a biquadratic number field $L'$ containing $\zeta_3$, there exist infinitely many elliptic curves $E_{a,b}$ and $\widehat{E}_{a,b}$ having
$$\dim_{\F_3} {\rm Sel}^{\widehat{\psi}_{a,b}}(\widehat{E}_{a,b}/K)-2 \le \dim_{\F_3} Cl(L')[3] \le \dim_{\F_3} {\rm Sel}^{{\psi}_{a,b}}({E}_{a,b}/K) \le \dim_{\F_3} {\rm Sel}^{3}({E}_{a,b}/K).$$
\end{corollary}
\begin{proof}
Note that $L' \cong \Q(\sqrt{-3}, \sqrt{a'})$ for some square-free integer $a'$. Since $-3a' \in L'$, we may assume $3 \nmid a'$. Set 
$T_{a'} := \big\{ \ell \in \Sigma_\Q \ : \ \ell \nmid 3a', \ \ell \equiv 1 \pmod{3} \ \text{  and  } \ \big(\frac{a'}{\ell} \big) =-1 \big\}$.
It is easy to see that $T_{a'}$ is an infinite set. Take $a := 16 a'$ if $a' \equiv 1 \pmod{4}$ and $a:=a'$ otherwise.
 We can verify that for such $a$ and for every $\ell \in T_{a'}$, the sets $S_1$ and $S_2$ in the Definition \ref{defofS1S2S3} are both empty.
 
For a fixed  $a$ as above,
consider the elliptic curves $E_{a,\ell_i}$ with $\ell_i \in T_{a'}$.  We claim that  $\{E_{a,\ell_i} \mid \ell_i \in T_{a'}\}$ is an infinite set of non-isomorphic elliptic curves over $K$. To see this, we compute the $j$-invariant  $j(E_{a,\ell_i})=-\frac{2^8a(a+6\ell_i)^3}{\ell_i^3(4a+27\ell_i)}$ and observe that there are at most four choices of $\ell_i$ such that $j(E_{a,\ell_i})$ are equal. 
For each of $\{E_{a,\ell_i} \mid \ell_i \in T_{a'}\}$, we get that $S_1 \cup S_2 = \emptyset$ and the result is immediate from Theorem \ref{containment2}.
\end{proof}

\newpage

\section*{Examples}
 Table \ref{tab:type2examples} contains the bounds on $\dim_{\F_3}{\rm Sel}^{\psi_{a,b}}(E_{a,b}/K)$ and $\dim_{\F_3}{\rm Sel}^{3}(E_{a,b}/K)$ of elliptic curves $E_{a,b}$ computed using SageMath/Magma. We define $s_l^\psi$ and $s_u^\psi$ to be the lower and upper bounds on the ${\rm Sel}^{\psi_{a,b}}(E_{a,b}/K)$, respectively, and $s_l^3$ and $s_u^3$ to be the lower and upper bounds on the ${\rm Sel}^3(E_{a,b}/K)$, respectively, following Theorem \ref{containment2}, Corollary \ref{type2sel3} and Theorem \ref{cortotype2}. 
 Also, $r$ denotes $\text{rk } E_{a,b}(K)$, which is computed using Magma under the assumption of \emph{GRH}. In the cases when Magma can not compute $r$ precisely, $[-,-]$ will denote the bounds on $r$ given by the software. 

\begin{table}[ht]
\centering
\renewcommand{\arraystretch}{0.95}
\begin{tabular}{ |c||c|c|c|c|c|c|c|c|c|c|}
\hline
\multirow{2}{*}{$a, b$} & \multirow{2}{*}{$S_1$} & \multirow{2}{*}{$S_2$} & \multirow{2}{*}{$S_3$} & \multirow{2}{*}{$h^3_{S_{1,2}(L)}$} & \multirow{2}{*}{$h^3_{S_{1,3}(L)}$} & \multirow{2}{*}{$r$} &   \multirow{2}{*}{$s_l^{\psi_{a,b}}$} & \multirow{2}{*}{$s_u^{\psi_{a,b}}$} & \multirow{2}{*}{$s^3_l$} & \multirow{2}{*}{$s^3_u$}\\
& & & & & & & & & & \\
\hline\hline
$79, 131$ & $\emptyset$ & $\{131\OO_K\}$ & $\emptyset$ & $2$ & $2$ & $[0,2]$ &  $2$ & $4$ & $2$ & $10$ \\
\hline
$137, 127$ & $\{2\OO_K\}$ & $\emptyset$ & $\{41\OO_K\}$ & $1$ & $1$ & $[1,2]$ & $1$ & $5$ & $1$ & $12$ \\
\hline
$137, 137$ & $\{2\OO_K\}$ & $\emptyset$ & $\emptyset$ & $1$ & $1$ & $1$ & $1$ & $4$ & $1$ & $9$ \\
\hline
$137, 143$ & $\{2\OO_K\}$ & $\{11\OO_K\}$ & $\{4409\OO_K\}$ & $1$ & $1$ & $1$ & $1$ & $6$ & $1$ & $13$ \\
\hline
$142, 83$ & $\emptyset$ & $\{83\OO_K\}$ & $\{53\OO_K\}$ & $2$ & $2$ & $3$ &  $2$ & $5$ & $3$ & $11$ \\
\hline
$142, 194$ & $\emptyset$ & $\emptyset$ & $\{2903\OO_K\}$ & $2$ & $2$ &  $[0,2]$ & $2$ & $4$ & $2$ & $8$ \\
\hline
$1714, 20$ & $\emptyset$ & $\{5\OO_K\}$ & $\emptyset$ & $3$ & $3$ & $2$ &  $3$ & $5$ & $3$ & $12$ \\
\hline
$2230, 48$ & $\emptyset$ & $\{\p\}$ & $\{1277\OO_K\}$ & ${2}$ & $3$ & $1$  & $2$ & $5$ & $2$ & $11$ \\
\hline
$2230, 533$ & $\emptyset$ & $\{41\OO_K\}$ & $\emptyset$ & $3$ & $3$ & $2$ &  $3$ & $5$ & $3$ & $12$ \\
\hline
$2263, 587$ & $\emptyset$ & $\{587\OO_K\}$ & $\emptyset$ & $3$ & $3$ & $2$  & $3$ & $5$ & $3$ & $12$ \\
\hline
$2659, 29$ & $\emptyset$ & $\{29\OO_K\}$ & $\emptyset$ & $3$ & $3$ & $2$ & $3$ & $5$ & $3$ & $12$ \\
\hline
$2659, 98$ & $\emptyset$ & $\emptyset$ & $\{29\OO_K\}$ & $3$ & $3$ & $2$ & $3$ & $5$ & $3$ & $10$ \\
\hline
$3023, 273$ & $\emptyset$ & $\emptyset$ & $\{19463\OO_K\}$ & $3$ & $3$ & $2$ & $3$ & $5$ & $3$ & $10$ \\
\hline
$3023, 678$ & $\emptyset$ & $\{113\OO_K\}$ & $\emptyset$ & $3$ & $3$ & $0$  & $3$ & $5$ & $3$ & $12$ \\
\hline
$79, 193$ & $\emptyset$ & $\emptyset$ & $\emptyset$ & $2$ & $2$ & ${[2,3]}$ &  $2$ & $3$ & $2$ & $7$ \\
\hline
$142, 12$ & $\emptyset$ & ${\{\p\}}$ & $\emptyset$ & ${1}$ & $2$ & $0$ &  $1$ & $3$ & $1$ & $7$  \\
\hline
$142, 67$ & $\emptyset$ & $\emptyset$ & $\emptyset$ & $2$ & $2$ & $1$  & $2$ & $3$ & $2$ & $7$\\
\hline
$223, 3$ & $\emptyset$ & $\{\p\}$ & $\emptyset$ & ${1}$ & $2$ & $0$  & $1$ & $3$ & $1$ & $7$ \\
\hline
$223, 63$ & $\emptyset$ & ${\{\p\}}$ & $\emptyset$ & ${1}$ & $2$ & $0$ & $1$ & $3$ & $1$ & $7$ \\
\hline
$2659, 24$ & $\emptyset$ & ${\{\p\}}$ & $\emptyset$ & ${2}$ & $3$ & $0$ &  $2$ & $4$ & $2$ & $9$ \\
\hline
$2659, 39$ & $\emptyset$ & ${\{\p\}}$ & $\emptyset$ & ${2}$ & $3$ & $0$ &  $2$ & $4$ & $2$ & $9$ \\
\hline
$3023, 18$ & $\emptyset$ & $\emptyset$ & $\emptyset$ & $3$ & $3$ & ${[1,3]}$ & $3$ & $4$ & $3$ & $9$ \\
\hline
$3391, 12$ & $\emptyset$ & ${\{\p\}}$ & $\emptyset$ & ${2}$ & $3$ & $0$ &  $2$ & $4$ & $2$ & $9$ \\
\hline
$3391, 67$ & $\emptyset$ & $\emptyset$ & $\emptyset$ & $3$ & $3$ & ${[0,1]}$ & $3$ & $4$ & $3$ & $9$ \\
\hline
$3667, 3$ & $\emptyset$ & ${\{\p\}}$ & $\emptyset$ & ${2}$ & $3$ & $0$ &  $2$ & $4$ & $2$ & $9$ \\
\hline
$3667, 63$ & $\emptyset$ & ${\{\p\}}$ & $\emptyset$ & ${2}$ & $3$ & ${[0,2]}$ &  $2$ & $4$ & $2$ & $9$ \\
\hline
$4094, 23$ & $\emptyset$ & $\emptyset$ & $\emptyset$ & $3$ & $3$ & ${[0,1]}$ &  $3$ & $4$ & $3$ & $9$ \\
\hline
$4094, 74$ & $\emptyset$ & $\emptyset$ & $\emptyset$ & $3$ & $3$ & $[0,1]$  & $3$ & $4$ & $3$ & $9$ \\
\hline
$4151, 14$ & $\emptyset$ & $\emptyset$ & $\emptyset$ & $3$ & $3$ & ${[0,3]}$ &  $3$ & $4$ & $3$ & $9$ \\
\hline
$4279, 39$ & $\emptyset$ & ${\{\p\}}$ & $\emptyset$ & ${2}$ & $3$ & ${[0,2]}$ & $2$ & $4$ & $2$ & $9$ \\
\hline
$4279, 66$ & $\emptyset$ & ${\{\p\}}$ & $\emptyset$ & ${2}$ & $3$ & $0$ & $2$ & $4$ & $2$ & $9$ \\
\hline
$4778, 48$ & $\emptyset$ & $\emptyset$ & $\emptyset$ & $3$ & $3$ & ${[0,1]}$ & $3$ & $4$ & $3$ & $9$ \\
\hline
$4910, 14$ & $\emptyset$ & $\emptyset$ & $\emptyset$ & $3$ & $3$ & ${[0,1]}$ & $3$ & $4$ & $3$ & $9$ \\
\hline
$4910, 38$ & $\emptyset$ & $\emptyset$ & $\emptyset$ & $3$ & $3$ & ${[0,1]}$ & $3$ & $4$ & $3$ & $9$ \\
\hline
$43063, 7$ & $\emptyset$ & $\emptyset$ & $\emptyset$ & $4$ & $4$ & ${[0,1]}$ & $3$ & $4$ & $3$ & $9$ \\
\hline
$43063, 96$ & $\emptyset$ & ${\{\p\}}$ & $\emptyset$ & ${3}$ & $4$ & $0$ & $3$ & $5$ & $3$ & $11$ \\
\hline
$51694, 13$ & $\emptyset$ & $\emptyset$ & $\emptyset$ & $4$ & $4$ & ${[1,3]}$ & $4$ & $5$ & $4$ & $11$ \\
\hline
$51694, 21$ & $\emptyset$ & ${\{\p\}}$ & $\emptyset$ & ${3}$ & $4$ & $0$ & $3$ & $5$ & $3$ & $11$ \\
\hline
$53507, 18$ & $\emptyset$ & $\emptyset$ & $\emptyset$ & $4$ & $4$ & ${1}$ & $4$ & $5$ & $4$ & $11$ \\
\hline
$53507, 74$ & $\emptyset$ & $\emptyset$ & $\emptyset$ & $4$ & $4$ & ${[0,1]}$ & $4$ & $5$ & $4$ & $11$ \\
\hline
$529987, 32$ & $\emptyset$ & $\emptyset$ & $\emptyset$ & $4$ & $4$ & ${[0,1]}$ & $4$ & $5$ & $4$ & $11$ \\
\hline
$529987, 108$ & $\emptyset$ & ${\{\p\}}$ & $\emptyset$ & ${3}$ & $4$ & $0$ & $3$ & $5$ & $3$ & $11$ \\
\hline
\end{tabular}
\caption{Bounds on the $\psi_{a,b}$- and $3$-Selmer groups of the curves $E_{a,b}:y^2=x^3+a(x-b)^2$}
\label{tab:type2examples}
\end{table}

\end{document}